\documentclass[preprint,1p,12pt]{elsarticle}
\usepackage[T1]{fontenc}
\usepackage[latin9]{inputenc}
\usepackage{amsthm}
\usepackage{amsmath}
\usepackage{amssymb}
\usepackage{graphicx}
\usepackage{setspace}
\usepackage{esint}
\usepackage{url}
\usepackage{comment}
\usepackage{color}

\makeatletter

\providecommand{\tabularnewline}{\\}
\newcommand{\lyxdot}{.}

\newcommand{\wasblue}[1] {#1}
\newcommand{\wasred}[1] {#1}

\newcommand{\refone}[1]{{\color{black}#1}}
\newcommand{\reftwo}[1]{{\color{black}#1}}

\makeatother

\begin{document}
\title{On the performance of exponential integrators for problems in magnetohydrodynamics}
\author[uibk]{Lukas Einkemmer\corref{cor1}} 	
\ead{lukas.einkemmer@uibk.ac.at}
\author[ucm]{Mayya Tokman}
\ead{mtokman@ucmerced.edu}
\author[llnl]{John Loffeld}
\ead{loffeld1@llnl.gov}
\address[uibk]{University of Innsbruck, \wasblue{Innsbruck}, Austria}
\cortext[cor1]{Corresponding author}
\address[ucm]{University of California, Merced, United States of America}
\address[llnl]{Lawrence Livermore National Laboratory, Livermore, United States of America}

\begin{abstract}
Exponential integrators have been introduced as an efficient alternative to explicit and implicit methods for integrating large stiff systems of differential equations.  Over the past decades these methods have been studied theoretically and their performance was evaluated using a range of test problems.  While the results of these investigations showed that exponential integrators can provide significant computational savings, the research on validating this hypothesis for large scale systems and understanding what classes of problems can particularly benefit from the use of the new techniques is in its initial stages.  Resistive magnetohydrodynamic (MHD) modeling is widely used in studying large scale behavior of laboratory and astrophysical plasmas. In many problems numerical solution of MHD equations is a challenging task due to the temporal stiffness of this system in the parameter regimes of interest. In this paper we evaluate the performance of exponential integrators on large MHD problems and compare them to a state-of-the-art implicit time integrator.  Both the variable and constant time step exponential methods of EpiRK-type are used to simulate magnetic reconnection and the Kelvin--Helmholtz instability in plasma. Performance of these methods, which are part of the EPIC software package, is compared to the variable time step variable order BDF scheme included in the CVODE (part of SUNDIALS) library.  We study performance of the methods on parallel architectures and with respect to magnitudes of important parameters such as Reynolds, Lundquist, and Prandtl numbers. We find that the exponential integrators provide superior or equal performance in most circumstances and conclude that further development of exponential methods for MHD problems is warranted and can lead to significant computational advantages for large scale stiff systems of differential equations such as MHD.  
\end{abstract}
\maketitle

\section{Introduction}
The focus of this article lies at the intersection of computational plasma physics and numerical analysis.  The system of resistive
magnetohydrodynamics (MHD) equations is used extensively in both laboratory and astrophysical plasma physics \cite{goedbloed2004}.  MHD simulations
are employed to study large scale plasma dynamics in applications ranging from \wasblue{the} design of the next generation of fusion devices to modeling
astrophysical jets (e.g. \cite{kuwabara2005, ebrahimi2015, zanni2007}).  Since experiments in this field usually impose a high cost and the ability to conduct observations can be very limited, computational modeling
in general, and numerical MHD in particular, have become very important components of research in plasma physics.  While MHD theory is widely used to study plasma dynamics, the complexity of the MHD equations makes their numerical solution a challenging task.  In particular,  the MHD system describes a wide range of temporal scales.  In applications such as large scale modeling of tokamaks or modeling of solar eruptions, the physical processes of interest evolve on the time scale of the shear Alfv\'en waves, which is much slower than the fastest mode in the system - the fast magnetosonic wave - and much faster than the slow magnetosonic wave and the resistive processes.  Such temporal stiffness coupled with the complexity of the equations poses a double-edged challenge to developing a time integrator for MHD.  From one perspective, the stiffness implies that explicit time integrators will suffer from prohibitively severe restrictions on the time step size due to stability considerations.  On the other hand the inherent three-dimensional nature of the processes  MHD describes and the complexity of spatial operators describing the dynamics makes evaluation of forcing terms computationally expensive. The latter difficulty served as one of the major reasons most large scale MHD codes came to rely extensively on explicit time integrators even for problems where resistivity plays a major role.  However, the temporal stiffness of the equations forced researchers working on these software packages to turn to at least some form of implicitness to overcome the small time steps sizes of explicit integrators \cite{jardin2012}.  Production codes used in fusion research and astrophysics \cite{NIMROD, M3D, FLASH} started embedding some form of implicitness into their time integration schemes.  The NIMROD code, for instance, employs an implicitly integrated term in the momentum equation derived by linearizing equations around a static magnetic field and using it to stabilize stiff wave components \cite{sovinec2010, NIMROD}.  In addition to this stabilizing term the overall time integration strategy also uses operator splitting and a staggered temporal grid.  The choice of the implicit terms in such semi-implicit schemes is also influenced by the availability of efficient algebraic solvers used to approximate solutions to the associated linear systems.  The complexity of the resulting overall time integration method makes it very difficult to use rigorous analysis to study stability and convergence.  Such an approach also precludes or makes it challenging to derive error estimators necessary for employing adaptive time stepping strategies or using the time integration as part of the uncertainty quantification methods.  Other partially implicit approaches such as alternating direction implicit (ADI) methods \cite{lindemuth73, schnack_killeen80} or implicit-explicit techniques \cite{keppens99} are also limited as they try to maintain a balance between stabilizing the time integration while maintaining the accuracy of the solution \cite{chacon08}

The difficulties associated with partially implicit methods and the ever pressing need for more efficient and accurate algorithms for high-performance computing platforms gave rise to the efforts to develop fully implicit time integrators for resistive MHD \cite{reynolds2006, reynolds2012, chacon08, lutjens10}.  A typical approach to constructing a fully implicit integrator involves using an implicit scheme such as, for instance, BDF and coupling it with a Newton--Krylov iteration to solve the implicit equations within each time step.  In \cite{reynolds2006, reynolds2012}, for example, the resistive MHD equations are solved using a variable order variable time step BDF method from the CVODE portion \cite{cohen1996} of the SUNDIALS package \cite{SUNDIALS}. While implicit methods possess better stability properties compared to explicit techniques, they are also affected by the stiffness of the equations.  In an implicit integrator temporal stiffness manifests itself in the stiff linear systems that have to be solved as part of each Newton iteration.  A key to obtaining an efficient fully implicit Newton--Krylov solver is in developing effective scalable preconditioners for these stiff linear systems.  For a complex anisotropic system such as MHD construction of a sufficiently efficient preconditioner is simultaneously a central and a challenging task that requires a significant investment of time and effort.  If the parameter regime is changed or an extra term is added to the system, it is possible that a new preconditioner has to be developed. Often the structure of the equations is exploited to build a physics based preconditioner \cite{chacon08} which makes the overall integrator highly dependent on the equations, thus precluding one from using off-the-shelf time integration software packages.

In recent years exponential integrators have emerged as a
promising alternative to explicit and implicit methods for solving evolution
equations. A significant body of research has been accumulated in
which these methods are investigated both from a theoretical (see
e.g. \citep{hochbruck2010} for a review article) as well as from
a numerical point of view (see e.g. \citep{tokman2010} and \citep{loffeld2012}).
This research includes the construction of exponential integrators
that have been tailored to a given differential equation (see e.g.
\citep{hochbruck1999}) as well as such integrators that can be applied
to a class of problems (see e.g. \citep{hochbruck2005, tokman2006}).
There is strong evidence to suggest that exponential integrators can
offer computational advantages compared to standard integrators particularly
for large scale complex problems where no efficient preconditioner is available.
Exponential integration can also be used as part of an Implicit-Exponential (IMEXP) scheme
\cite{tokman14} for problems where an efficient
preconditioner has been constructed for some terms in the equations.

Exponential integrators were first shown to be a promising approach in the context of the MHD equations in \cite{tokman2001, tokman2002}. Since then
a number of more efficient and sophisticated exponential time integrators have been developed.
Exponential propagation iterative methods of Runge--Kutta type (EpiRK) were
first introduced in \cite{tokman2006} for general nonlinear systems of ordinary differential
equations.  The EpiRK framework \cite{tokman2006, tokman2010, tokman14} for constructing exponential time integrators has been
designed to reduce computational complexity in each time step and provide a flexible ansatz
that allows construction of particularly efficient exponential schemes \cite{tokman2010, tokman2011, rainwater14, tokman14}.
The most efficient EpiRK methods have been implemented as part of the Exponential Integrator
Collection (EPIC) package, a C++/MPI software developed for serial and parallel computational
platforms (available for download by request to mtokman@ucmerced.edu).  Performance of EPIC integrators has been tested on a suit of numerical examples and it was shown that the EpiRK schemes can achieve superior or equal performance
compared to the BDF integrator found in the CVODE library \cite{loffeld2012, loffeld14}.

However, up to now most numerical tests of exponential methods found in the literature considered
problems of small to medium complexity (see, for example, \citep{loffeld2012},
\citep{klein2011}, or \citep{einkemmer2013}) where in many instances
efficient preconditioners can be constructed. The goal of this paper is to apply the EpiRK integrators found in the
EPIC library to the problem of solving the resistive magnetohydrodynamics
(MHD) equations. General preconditioners for such systems have, for example,
been investigated in \citep{reynolds2012}
and it was found that constructing an efficient preconditioner for such system is challenging and the
performance gain from standard preconditioning techniques for
several dimensional examples is rather modest under most circumstances. Clearly a problem specific
preconditioner or a tailor-made time integrator for a specific system of equations will be the most efficient
approach in the majority of cases, but developing such a scheme requires extensive effort and time. The goal
of this paper is to investigate what computational savings can be achieved with a general exponential integrator
compared to an unpreconditioned implicit scheme. Thus we will perform a comparison of EPIC integrators with the BDF methods from
the CVODE library without a preconditioner.  We aim to demonstrate that exponential integration is a viable technique to use
in integrating MHD equations.

The paper is organized as follows.  In section \ref{sec:Exponential-integrators} we provide an introduction
to exponential integrators in general as well as to the specific EpiRK
methods we use in the numerical examples and their implementation in EPIC. In section \ref{sec:magnetohydrodynami}
we describe the MHD equations which are solved numerically in section
\ref{sec:The-reconnection-problem} for a reconnection problem, and
in section \ref{sec:The-Kelvin-Helmholtz-instability} for the Kelvin--Helmholtz
instability. Implementation details are discussed in section \ref{sec:Implementation}.
Finally, our conclusions are presented in section \ref{sec:Conclusion}.

\section{Exponential integrators\label{sec:Exponential-integrators}}

\wasblue{Let us consider the following initial value problem
\begin{align}
&y^{\prime}  =F(y)\label{eq:ivp}, \\
&y(t_0)  =y_{0}\nonumber
\end{align}
which is assumed to be large ($y, F(y) \in \mathbf{R}^N$ with $N \reftwo{\gg} 1$) and stiff.  Let us define by $y(t_n)$ and $y(t_n+h)$ solutions to \eqref{eq:ivp} at times $t_n$ and $t_n+h$ respectively.  Assuming that the first order Taylor expansion of the solution $y(t_n+h)$ around solution $y(t_n)$ exists we can 
write $y(t_n+h)$ in an integral form as follows. Let us denote the remainder function of the first order Taylor expansion of $F(y)$ at \refone{$\overline{y}$ as follows 
\begin{equation}
    R(y,\overline{y}) = F(y) - F(\overline{y}) - F'(\overline{y}) (y-\overline{y}).
\end{equation} }
Then \eqref{eq:ivp} can be written as 
\begin{equation}
    y' = F(y) = F(y(t_n)) + F'(y(t_n))(y-y(t_n)) + \refone{R(y,y(t_n))}. 
    \label{eq:taylor}
\end{equation}
To simplify the notation we define \refone{$J(t_n)=F'(y(t_n))$} and \refone{$F(t_n) = F(y(t_n))$}. Multiplying equation \eqref{eq:taylor} by a factor \refone{$e^{tJ(t_n)}$} and integrating the resulting formula from $t_n$ to $t_n+h$, we obtain the variation-of-constants formula
\begin{equation}
    \refone{
        \begin{aligned}
            y(t_n+h)&=y(t_n)+\left(e^{hJ(t_n)}-I\right)(hJ(t_n))^{-1}hF(t_n)\\
            &\quad+\int_{t_n}^{t_n+h}e^{(t_n+h-t)J(t_n)}R(y(t),y(t_n))\,dt.
        \end{aligned} 
    }
        \label{eq:intivp1}
\end{equation}
Changing the integration variable in \eqref{eq:intivp1} to $s$ with $t = t_n + sh$, $0 < s < 1$ we obtain
\begin{equation}
    \refone{
        \begin{aligned}
            y(t_n+h) &=y(t_n)+\left(e^{hJ(t_n)}-I\right)(hJ(t_n))^{-1}hF(t_n)\\
            &\quad+\int_{0}^{1}e^{h(1-s)J(t_n)}hR(y(t_n+sh),y(t_n))\,ds.
        \end{aligned} 
    } 
        \label{eq:intivp}
\end{equation}
If $t_n$ and $t_n+h$ represent discretization nodes on the time interval $[t_0, t_{end}] $ over which we want to approximate $y(t)$ ($t \in [t_0,t_{end}]$) we can 
use the integral representation of the solution $y(t_n+h)$ in \eqref{eq:intivp} to construct an exponential integrator.}
\refone{In the following we will denote the numerical approximation at time $t_n$ by $y_n$, i.e.~$y_n \approx y(t_n)$ and our goal is thus to compute $y_{n+1} \approx y(t_{n+1})$ given $y_n$.}

\wasblue{
\refone{To do that equation (\ref{eq:intivp}) has to be approximated in a suitable manner. This task is done in two steps. First, since we do not know $y(t_n)$, $F(t_n)$, and $J(t_n)$ we replace them by the numerical approximations $y_n$, $F_n=F(y_n)$, and $J_n=J(y_n)$. Then we use a polynomial approximation of the function $R(\cdot,y_n)$} on some nodes $s_k$. The integral in \eqref{eq:intivp} can then be written as a linear combination of terms of type $\varphi_k(s_khJ_n)v_k$, where the entire functions $\varphi_k(z)$ are defined by
\begin{align}
&\varphi_0(z)  =e^{z} \nonumber \\
&\varphi_k(z)  =\int_0^1 e^{(1-s)z} \frac{s^{k-1}}{(k-1)!}ds.  \label{eq:varphi}
\end{align}
Given definitions \eqref{eq:varphi} the second term on the right-hand-side of \eqref{eq:intivp} can be written as \refone{$\varphi_1(hJ_n)h F(t_n)$}.}

\wasblue{Note that all terms of type $\varphi_k(chJ)v$, where $c$ and $h$ are constants, $J$ is a matrix and $v$ is a vector, are products of matrix functions and vectors.  
Since $J$ is the Jacobian of the large stiff system \eqref{eq:ivp}, it is unlikely that diagonalizing the matrix $J$ is a computationally cheap, or even a feasible task. Thus the terms  $\varphi_k(chJ)v$ have to be approximated in a different way.  Computing these terms constitutes the main computational cost of an exponential method. Therefore choosing the quadrature for the integral in \eqref{eq:intivp} in order to minimize the required number of these calculations as well as selecting an appropriate algorithm to evaluate terms of type $\varphi_k(chJ)v$ is essential to the construction of an efficient exponential integrator. Note that construction of a particular integrator involves choosing the parameter $c$. Small values of $c$ can reduce computational cost as explained below. To explicitly make this point we leave the argument of the $\varphi_k$ function as $chJ$ rather than use a simplified notation with a single matrix variable. }

A number of methods have been proposed in the literature to
\wasred{approximate the exponential-like matrix functions $\varphi_k(chJ)$ or their product with a given vector}
%approximate either a general or a specific (e.g. rational or exponential) function of a matrix or its product with a given vector
(e.g. \wasblue{\cite{niesen2012, almohy2011,kassamtrefethen,leja,MN04,EH06}} or see the review in \cite{hochbruck2010}).  However, many of these techniques are applicable only to small matrices $J$ \cite{molervanloan, sidje}. In this paper we are concerned with large scale problems such as MHD, where
the dimensionality of the system \eqref{eq:ivp} is high (e.g. number of unknowns in \wasblue{the} MHD equations multiplied by the number of grid points in two- or three-dimensions),
and thus standard techniques for approximating a product of a matrix function $\varphi_k$ with a vector are either computationally infeasible or prohibitively expensive.

\reftwo{In this paper we use an approach} based on the Krylov projection algorithm \cite{arnoldi,vandervorst}.
Krylov methods have been extensively used to approximate the inverse of a matrix \cite{saad} and in 1987 were employed by Van der Vorst to approximate
a general function of a matrix and a vector product \cite{vandervorst}.  The algorithm works
by iteratively computing a projection of $\varphi_k(chJ)v$ onto a Krylov subspace
\[
    \refone{K_{m}(J,v)=\text{span}\left\{ v,Jv,\dots,J^{m-1}v\right\}}
\]
using the Arnoldi algorithm \cite{arnoldi}. \refone{Note that the Krylov subspace does not depend on $c$ and $h$ but since we are computing $\varphi_k(chJ)v$ the number of
Krylov vectors computed $m$ depends on the parameters $c$ and $h$. In fact, the value of} $m$ is determined dynamically by estimating the residuals in
the course of the iteration \cite{saadresid, hochbruckkrylov}.  For each $m$ a matrix $V_m$ is
formed that consists of the orthonormal vectors that form the basis of $K_m(chJ, v)$.
A side product of the Arnoldi iteration is the $m \times m$ matrix $H_m = chV_{m}^{T}JV_{m}$ which
represents the projection of $chJ$ onto the Krylov subspace.  A general product of the
matrix function $f$ with a vector is then approximated as
\[
f(chJ)b\approx V_{m}\left(V_{m}^{T}\varphi_{k}(chJ)V_{m}\right)V_{m}^{T}v \approx V_{m}\varphi_{k}(H_{m})V_{m}^{T}v.
\]
\wasblue{$H_{m}$ is a matrix of dimension $m\times m$ and a small $m$ is usually sufficient to achieve the desired accuracy. Then 
the application of $\varphi_{k}(H_{m})$ to a vector can be cheaply
computed by using any of the standard methods (such as Pad\'e approximation
or methods based on polynomial interpolation \cite{molervanloan, hochbruck2010}). The residuals are used to determine the stopping value of $m$ that meets the provided tolerance.  Note that the Arnoldi iteration portion of the algorithm is a scale invariant procedure.  Thus several terms $\varphi_k(c_ihJ)v$ ($i = 1,2,...,K$) can be computed using a single Arnoldi iteration, i.e. using the same Krylov basis.  To meet the desired accuracy requirements we just have to ensure that the dimensionality of the basis $m$ is chosen with respect to computing $\varphi_k(c_ihJ)v$ with the largest $c_i$. This can be accomplished using appropriate residuals as described for example in \cite{tokman2010}.}

Another possibility is to use a direct polynomial approximation of $\varphi_{k}$. To that end we have to chose $m$
interpolation nodes on some compact set $K\subset\mathbb{C}$. Chebyshev
nodes are an obvious choice, however, they suffer from the disadvantage
that to compute the interpolation for $m+1$ points we have to reevaluate
all the matrix-vector products already computed. To remedy this, Leja
points have been proposed which share many of the favorable properties
of the Chebyshev nodes but can be generated in sequence \wasred{\cite{leja,martinez2009,caliariostermann,leja_relpm,CKOR14,CKOR16}}.
Another advantage of interpolation at Leja points is its modest memory requirements. This makes
it an attractive alternative for computer systems where memory is
limited (such as graphic processing units, see \citep{einkemmer2013}).
The main disadvantage of the method, however, is that (at least) some
approximation of the spectrum of $J$ has to be available. 

\wasred{We also mention another algorithm for approximating products of type $\varphi_k(chJ)v$ -- the Taylor polynomial approximation-based 
method introduced in \cite{almohy2011} \reftwo{(in fact, this is a direct polynomial interpolation method, as described above, where all interpolation nodes coincide in a single point)}. This technique combines ideas of scaling-and-squaring, Taylor expansion and 
certain error approximations based on the 1-norm of the matrix to estimate $\varphi_k(chJ)v$ or sequences of $\varphi_k(\tau_i chJ)$ at equally 
spaced intervals $\tau_i$, $\tau_i = i\tau$. As indicated in the experiments conducted in \cite{almohy2011} (e.g. Experiment 4), the adaptive Krylov algorithm of \cite{niesen2012} can outperform the \reftwo{Taylor method}. Note that not all numerical experiments in \cite{almohy2011}, particularly Experiment 9 which involves the 
largest matrix (dimension 250,000), include comparisons with the adaptive Krylov method of Niesen and Wright \cite{niesen2012}. 
In our numerical tests we found that for the large scale matrices we encountered in MHD and other test problems the adaptive Krylov method was more efficient than the Taylor-based algorithm.  We have not encountered the instability in the adaptive Krylov algorithm found in Experiment 10 in \cite{almohy2011} because such instability can only occur with a Krylov-based method if $p$ and $c$ are increased without bound in the linear combination $\sum_{k=1}^p\varphi_k(chJ)v_k$ (typically $p > 10$). This case is not relevant to most of the high-order exponential integrators proposed in the literature including EpiRK methods which use linear combination of only a few $\varphi_k$ ($p \le 5$) with appropriately (and adaptively) scaled $c$ and $h$. Additionally, as described above, the ability of the adaptive Krylov methods to yield an approximation for any discretization $\tau_i$ (not only the equally spaced $\tau_i$'s) allows taking advantage of the EpiRK framework to reduce the number of Krylov projections and to make each projection more efficient by optimizing the coefficients $c$. Unlike the Taylor-based method, the adaptive Krylov algorithm does not require any estimates of the matrix norm of $chJ$ which reduces computational cost and is preferable for matrix-free implementations. \reftwo{We note, however, that for problems where the matrix is stored in memory Taylor methods or methods based on Leja interpolation can be beneficial (for example, a backward error analysis is available for these methods \cite{almohy2011,CKOR16}).}  Since approximating the action of matrix functions is an active area of research, it is likely that more efficient algorithms, including Taylor-based methods, will emerge in the future that allow us to further reduce the computational cost of exponential integrators. At present we focus on the EpiRK methods based on an adaptive Krylov approach.}

Since evaluations of $\varphi_k(chJ)v$ constitute the most computationally expensive part of the exponential
integrator, special care has to be taken to ensure that each time step requires as few of these evaluations as possible
and the most efficient algorithm is chosen for these approximations.  The EpiRK framework \cite{tokman2006, tokman2011}
is designed to construct integrators with a reduced number of $\varphi_k(chJ)v$ evaluations and with $\varphi_k$ functions
chosen to speed up the Krylov iteration. \wasblue{In addition,  the EpiRK framework allows optimizing the coefficients $c$ to gain
more efficiency. Recall that the Krylov iteration converges faster if the spectrum of the matrix $chJ$ is more clustered. Minimizing the coefficients $c$
increases clustering of the spectrum of $chJ$ and consequently reduces the dimensionality $m$ of the Krylov subspace for a given tolerance.} 

A particularly efficient class of EpiRK methods \cite{tokman2010} also takes advantage
of an adaptive version of the Krylov iteration introduced in \cite{niesen2012}.  The adaptive Krylov based EpiRK schemes link
the construction of a time integrator with the choice of the method to approximate $\varphi_k(chJ)v$ to build a more efficient
overall time stepping scheme in the following way.

The computational cost of the Krylov algorithm to approximate terms of type $\varphi_k(c hJ)v$ scales as $O(m^2)$ with the size of the
Krylov basis $m$. Clearly, $m$ depends on the eigenvalues of $J$, \wasblue{the} vector $v$ as well as the values of $c$ and $h$.
As $c$ or $h$ is increased, so is the size of the basis.  The EpiRK framework allows one to reduce the coefficients $c$ and choose
the better suited functions $\varphi_k$ in the derivation of a time integrator.  The adaptive Krylov technique allows replacing each computation
of $\varphi_k(c hJ)v$ with several evaluations of $\varphi_k(\tau_i c hJ)v$ where $0 < \tau_i < 1$ for $i = 0, 1, 2, ... \wasblue{,P}$.   Since the complexity of each of
these evaluations now scales as $O(m_i^2)$, in most cases increased efficiency is observed as $O(m^2)$ can be \wasblue{larger} than $O(m_1^2)+O(m_2^2)+...+O(m_P^2)$.

A general EpiRK method can be written as
\begin{align*}
    &Y_{i}  =y_{n}+a_{i1}\psi_{i1}(g_{i1}hJ_{i1})hF_{n}+\sum_{j=2}^{i-1}a_{ij}\psi_{ij}(g_{ij}hJ_{ij})h\Delta^{(j-1)}\refone{R_n},  \quad i = 1,..,s-1\\
    &y_{n+1}  =y_{n}+b_{1}\psi_{s1}(g_{s1}hJ_{s1})hF_{n}+\sum_{j=2}^{s}b_{j}\psi_{sj}(g_{sj}hJ_{sj})h\Delta^{(j-1)}\refone{R_n}
\end{align*}
where \refone{$R_n(\cdot) = R(\cdot,y_n)$ and} the divided differences \refone{$\Delta^{(j-1)}R_n$}  are computed by using the nodes $y_{n},Y_{1},\dots,Y_{s-1}$ (note that since \refone{$R_n(y_n) = 0$} the node $y_n$ does not actually appear in the divided differences).  Different choices
of functions $\psi_{ij}$ and matrices $J_{ij}$ result in different classes of EpiRK methods such as unpartitioned, partitioned or hybrid
exponential or implicit-exponential (IMEXP) integrators \cite{tokman14}.  In this paper we employ the schemes with the choice of $J_{ij}=J_n$ and
$\psi_{ij}$ set as linear combinations of $\varphi_k$ functions
\[
\psi_{ij}(z)=\sum_{k=1}^{s}p_{ijk}\varphi_{k}(z).
\]

As in the Runge--Kutta case, to construct an efficient integrator the
coefficients $a_{ij},g_{ij},b_{j}$ and $p_{ijk}$ have to be chosen
subject to the appropriate order conditions.   In this paper, we use
the variable time step fifth order method EpiRK5P1
\begin{align}
&Y_1 = y_n + a_{11}\varphi_{1}(g_{11}hJ_n)hF_n\label{eqn:epirk5p1} \\
&Y_2 = y_n + a_{21}\varphi_{1}(g_{21}hJ_n)hF_n + a_{22}\varphi_{1}(g_{22}hJ_n)hR(Y_1) \nonumber \\
&y_1 = y_n + b_1\varphi_1(g_{31}hJ_n)hF_n + b_2\varphi_1(g_{32}hJ_n)hR(Y_1) + b_3\varphi_3(g_{33}hJ_n)h(-2R(Y1)+R(Y2)),\nonumber
\end{align}
that has been derived in \citep{tokman2010} (see Table \ref{tbl:epirk5p1} for the coefficients; this
method is also employed in the performance comparison given in
\citep{loffeld2012}) and the constant time step fourth order EpiRK4 method
proposed in \cite{rainwater2016} given by
\begin{equation}
\label{eqn:epirk4}
  \begin{aligned}
&Y_{1}= y_n+\tfrac{1}{2}\varphi_1(\tfrac{1}{2}h_nJ_n)h_nF_n\\
&Y_{2}= y_n+\tfrac{2}{3}\varphi_1(\tfrac{2}{3}h_nJ_n)h_nF_n \\
&y_{n+1} = y_n + \varphi_1(h_nJ_n)h_nF_n+\left(32\varphi_3(h_nJ_n)-144\varphi_4(h_nJ_n)\right)h_nR(Y_{1})\\
&\quad \quad \quad +\left(-\tfrac{27}{2}\varphi_3(h_nJ_n)+81\varphi_4(h_nJ_n)\right)h_n R(Y_{2})
  \end{aligned}
  \end{equation}

\begin{table}[htp]
\caption{Coefficients of the fifth order method EPIRK5P1.}
\footnotesize
\begin{tabular}{|l|}
\hline \\
$\begin{bmatrix} a_{11} & & \\ a_{21} & a_{22} & \\ b_1 & b_2 & b_3   \end{bmatrix}$ =
$\begin{bmatrix}
 0.35129592695058193092 & &\\
 0.84405472011657126298 & 1.6905891609568963624 & \\
1.0 & 1.2727127317356892397 &   2.2714599265422622275
\end{bmatrix}$ \\[25pt]
%\hline
$\begin{bmatrix} g_{11} & & \\ g_{21} & g_{22} & \\ g_{31} & g_{32} & g_{33}   \end{bmatrix}$ =
% Coefficients g
$\begin{bmatrix}
 0.35129592695058193092& &\\
 0.84405472011657126298 & 0.5 & \\
1.0 & 0.71111095364366870359 &  0.62378111953371494809
\end{bmatrix}$  \\[25pt]
% Coefficients a, b
\hline
\end{tabular}
\label{tbl:epirk5p1}
\end{table}
Both of these EpiRK schemes utilize the adaptive Krylov
technique for the evaluations of $\varphi_k(chJ)v$ in order to improve performance.  It is important to note that while
both of these are three-stage methods, EpiRK4 requires only two Krylov projections to be executed
per time step since the scale invariance of the Arnoldi iteration allows \wasblue{the} simultaneous calculation of the
Krylov bases for $Y_1$ and $Y_2$ and one more Krylov projection to evaluate $y_{n+1}$ (e.g. see \cite{rainwater2016}
for details). EPIRK5P1 requires three Krylov projections to be performed each time step \cite{tokman2010}.
In addition EPIRK4 is a so-called stiffly accurate method while EPIRK5P1 is derived using classical order
conditions.   The advantages of the EPIRK5P1 method, however, include the higher order and the ability to
design embedded methods for the automatic time step control that do not require an additional Krylov projection.
The adaptive time step control version of the methods of type EPIRK4 is still under development.

\section{The resistive magnetohydrodynamics equations\label{sec:magnetohydrodynami}}

The most fundamental theoretical description of a classical plasma
comes from the kinetic equation. This so-called Vlasov equation (the
collisionless case) or Boltzmann equation (when collisions are of
physical significance) describes the time evolution of a particle-density
in the 3 + 3 dimensional phase space (the first three dimensions correspond
to the space dependence while the remaining correspond to the velocity
dependence of the particle-density). While a number of simulations
of different plasma phenomena have been conducted with this approach,
due to the high dimensionality of the phase space, a lower-dimensional
approximation is usually used to render such simulations feasible.
A similar discussion holds true for the gyrokinetic approximation
that can be employed in plasmas where a strong external magnetic field
along a given axis is present. In this case the phase space of the
Vlasov equation is reduced to $3+2$ dimensions by averaging over
the gyro motion. This procedure yields a good approximation under
the assumption of low frequency as compared to the cyclotron frequency.
For a more detailed discussion of gyrokinetic models, see e.g. \citep{frieman1982}
or \citep{fahey2004}.

However, in many applications (such as magnetic confinement fusion,
spheromak experiments, and astrophysical plasmas) the timescales of
interest are sufficiently long and/or the full three dimensional model
is necessary to describe the physical phenomena. In that case the
kinetic approach is usually not feasible (even on modern day supercomputers)
and thus further simplifications have to be introduced. For the magnetohydrodynamics
(MHD) equations, which we will describe in the remainder of this section,
the assumption is made that, to a good approximation, the distribution
in the velocity direction is Maxwellian; in other words, it is assumed that each sufficiently
small volume in the plasma is in thermodynamic equilibrium.  For a general overview of kinetic and MHD models (including a derivation
of the MHD equations from the Vlasov equation) see e.g. \citep{nicholson1983}.
Numerical computations in the context of the MHD model discussed in
this paper are performed, for example, in \citep{reynolds2006}, \citep{reynolds2010},
and \citep{reynolds2012}.

If the assumption is made that the plasma considered is in thermodynamic
equilibrium, the equations of motion in a three dimensional phase
space (a so-called fluid model) can be derived. These are given by
(in dimensionless units)
\begin{align}
\frac{\partial\rho}{\partial t}+\nabla\cdot(\rho\boldsymbol{v}) & =0\label{eq:continuity}\\
\rho\frac{\partial\boldsymbol{v}}{\partial t}+\rho\boldsymbol{v}\cdot\nabla\boldsymbol{v} & =\boldsymbol{J}\times\boldsymbol{B}-\nabla p\label{eq:momentum-equation}
\end{align}
which we refer to as the continuity and momentum equation respectively.
The equations are written in terms of the density $\rho$, the fluid
velocity $\boldsymbol{v}$, the magnetic field $\boldsymbol{B}$,
the electric current density $\boldsymbol{J}$, and the pressure
$p$. The dynamics is determined by the Lorentz force (the $\boldsymbol{J}\times\boldsymbol{B}$
term) and the pressure gradient force (the $\nabla p$ term). These
fluid equations have to be coupled to an appropriate model of the
electric field. Note, however, that if \wasblue{(the ideal)} Ohm's law is assumed
to hold, i.e.
\[
\boldsymbol{E}+\boldsymbol{v}\times\boldsymbol{B}=0,
\]
then the electric field $\boldsymbol{E}$ can be eliminated from Maxwell's
equation. This leaves us with the following system
\begin{align}
	&\wasblue{\frac{\partial\boldsymbol{B}}{\partial t}  =\nabla\times(\boldsymbol{v}\times\boldsymbol{B})}\label{eq:Maxwell}\\
&\boldsymbol{J} =\nabla\times\boldsymbol{B}\nonumber \\
&\nabla\cdot\boldsymbol{B}  =0.\nonumber
\end{align}
The first of these equations yields the time evolution of $\boldsymbol{B}$,
the second can be used to eliminate $\boldsymbol{J}$ from equation
(\ref{eq:momentum-equation}), and the third is the familiar solenoidal
constraint imposed on the magnetic field.

A commonly employed approach to close these equations (see e.g. \citep{reynolds2010})
is to supplement them with the following equation of state
\[
e=\frac{p}{\gamma-1}+\frac{\rho}{2}v^{2}+B^{2},
\]
where the time evolution of the energy density $e$ is given by
\begin{equation}
	\wasblue{\frac{\partial e}{\partial t}+\nabla\cdot\left((e+p+\tfrac{1}{2}B^{2})\boldsymbol{v}-\boldsymbol{B}(\boldsymbol{B}\cdot\boldsymbol{v})\right)=0.}\label{eq:energy}
\end{equation}

Collectively, the equations (\ref{eq:continuity})-(\ref{eq:energy})
in the variables $(\rho,\boldsymbol{v},\boldsymbol{B},e)$ yield a
first-order system of $8$ differential equations in $8$ variables
and are called the ideal magnetohydrodynamics equations (or the ideal
MHD equations).

For the purpose of performing the spatial discretization, these equations
are often cast into the so-called divergence form. Then, the equations
of motion read as (see e.g. \citep{reynolds2010})
\begin{align*}
	\frac{\partial U}{\partial t}+\nabla\cdot \wasred{F_h(U)} & =0
\end{align*}
with state vector
\[
U=\left[\begin{array}{c}
\rho\\
	\wasblue{\rho \boldsymbol{v}}\\
\boldsymbol{B}\\
e
\end{array}\right]
\]
and
\begin{equation}
	\wasred{F_h(U)}=\left[\begin{array}{c}
\rho\boldsymbol{v}\\
\rho\boldsymbol{v}\otimes\boldsymbol{v}+(p+\tfrac{1}{2}B^{2})I-\boldsymbol{B}\otimes\boldsymbol{B}\\
\boldsymbol{v}\otimes\boldsymbol{B}-\boldsymbol{B}\otimes\boldsymbol{v}\\
(e+p+\tfrac{1}{2}B^{2})\boldsymbol{v}-\boldsymbol{B}(\boldsymbol{B}\cdot\boldsymbol{v})
\end{array}\right],\label{eq:F}
\end{equation}
where we have denoted the tensor product by using the $\otimes$ symbol.

In this paper, as in \citep{reynolds2006} and \citep{reynolds2010},
we will consider a slightly more general class of equations which,
in addition to the dynamics discussed so far, includes dissipative
effects (due to particle collisions in the plasma). To that end the
hyperbolic flux vector \wasred{$F_h(U)$} is extended in \citep{reynolds2006}
by a diffusive part given by
\begin{equation}
F_{d}(U)=\left[\begin{array}{c}
0\\
Re^{-1}\boldsymbol{\tau}\\
S^{-1}\left(\eta\nabla B-\eta(\nabla B)^{\mathrm{T}}\right)\\
Re^{-1}\boldsymbol{\tau}\cdot\boldsymbol{v}+\frac{\gamma}{\gamma-1}Re^{-1}Pr^{-1}\nabla T+S^{-1}\left(\tfrac{1}{2}\nabla(\boldsymbol{B}\cdot\boldsymbol{B})-B(\nabla\boldsymbol{B})^{\mathrm{T}}\right)
\end{array}\right]\label{eq:Fd}
\end{equation}
with
\[
\boldsymbol{\tau}=\nabla\boldsymbol{v}+(\nabla\boldsymbol{v})^{T}-\tfrac{2}{3}\nabla\cdot\boldsymbol{v}I,
\]
where we have assumed a spatially homogeneous
viscosity $\mu$, resistivity $\eta$, and thermal conductivity $\kappa$.
The non-dimensional parameters governing the solution of this system are
the Reynolds number \wasblue{$Re=\rho_{0}v_{A}l_{0}/\mu_0$},
the Lundquist number \wasblue{$S=\mu_{0}v_{A}l_{0}/\eta_0$} and the Prandtl
number \wasblue{$Pr=c_{p}\mu_0/\kappa_0$} for a characteristic density $\rho_{0}$,
a characteristic length scale $l_{0}$, and the Alfven velocity $v_{A}$
(as usual the permeability of free space is denoted by $\mu_{0}$,
\wasblue{$\kappa_0=5/3$}, and $c_{p}$ is the specific heat of the fluid).

Note that in all the simulations we use dimensionless units. Thus, we can rewrite (\ref{eq:Fd}) more conveniently in terms
of the (dimensionless) viscosity $\mu=Re^{-1}$, the (dimensionless)
resistivity $\eta=S^{-1}$, and the (dimensionless) thermal conductivity
$\kappa=Pr^{-1}$. Then, the form of the resistive MHD equations used
for the spatial discretization is
\[
	\frac{\partial U}{\partial t}+\nabla\cdot \wasred{F_h(U)}=\nabla\cdot F_{d}(U)
\]
with \wasred{$F_h(U)$} given by equation (\ref{eq:F}) and where $F_{d}(U)$
is given by
\[
F_{d}(U)=\left[\begin{array}{c}
0\\
\mu\boldsymbol{\tau}\\
\eta\left(\nabla B-(\nabla B)^{\mathrm{T}}\right)\\
\mu\boldsymbol{\tau}\cdot\boldsymbol{v}+\frac{\gamma\mu\kappa}{\gamma-1}\nabla T+\eta\left(\tfrac{1}{2}\nabla(\boldsymbol{B}\cdot\boldsymbol{B})-B(\nabla\boldsymbol{B})^{\mathrm{T}}\right)
\end{array}\right].
\]
For the numerical simulations conducted in this paper, we employ the MHD code
developed in \cite{einkemmercppmhd}. This C++ code is based on the Fortran code developed in
\citep{reynolds2006} which has been used to conduct plasma physics simulations
(see e.g. \citep{reynolds2006}, \citep{reynolds2010}, and
\citep{reynolds2012}) as well as to construct more
efficient preconditioners in the context of implicit time integrators (see
\citep{reynolds2010}).
\wasred{The implementation of exponential integrators requires the evaluation of the action of the Jacobian applied to a vector. We compute this action by employing a simple forward difference stencil (this has the advantage that one function evaluation can be reused). Consequently, the present implementation is completely matrix-free. For more details we refer the reader to section \ref{sec:Implementation}.}

\wasred{With respect to space discretization,} our implementation assumes a
finite difference or finite volume method where a single value is stored in each
cell. Thus, in each cell we store the value of the density $\rho$, the fluid
velocity $\boldsymbol{v}$, the magnetic field $\boldsymbol{B}$, the energy $e$
(but not the pressure p). For the numerical simulations conducted in
this paper, we have implemented the 2.5-dimensional case (that is, the state
variables only depend on the $x$- and $y$-direction but both the velocity
$\boldsymbol{v}(x,y)$ and the magnetic field $\boldsymbol{B}(x,y)$ are
three-dimensional vectors) using a classic centered stencil for the divergence;
that is, for each vector $G$, corresponding to the flux of a given (scalar) state variable, we compute the following approximation of the divergence
\[ \nabla\cdot
	G(U)\approx\frac{G(U_{i+1,j})-G(U_{i-1,j})}{2h}+\frac{G(U_{i,j+1})-G(U_{i,j-1})}{2h},
\]
where $i$ and $j$ are the cell indices in the $x$- and $y$-directions,
respectively and $h$ is the cell size. To evaluate the spatial derivatives
present in the diffusion vector $F_{d}$, we also employ the classic centered
stencil. Thus, the code has order two accuracy in space.

\wasred{It is well known that for MHD problems preserving the divergence free property
of the magnetic field is important. The space discretization used in our
implementation ensures that if the solenoidal property $\nabla\cdot\boldsymbol{B}=0$
is satisfied for the initial value, then this is true for all later
times. In \citep{reynolds2010} it is shown that this property also holds true if we couple the
space discretization with an implicit method that is based on a matrix-free
inexact Newton\textendash Krylov algorithm.
}

\wasred{ We will now show that an equivalent result holds for the exponential integrators considered in this paper. That is, if the discrete selenoidal property 
\[ \frac{B_{i+1,j}(t)-B_{i-1,j}(t)}{2h} + \frac{B_{i,j+1}(t)-B_{i,j-1}(t)}{2h} = 0 \]
is satisfied for the initial value this is also true for all later times (independent of the step size and the tolerance specified for the computation of the $\varphi_k$ functions). 
Now, it is well known that $F(U)$ satisfies the discrete selenoidal condition if $U$ does (see, for example, \citep{reynolds2010}). In addition, since the Jacobian $J$ is computed using a finite difference approximation, the selenoidal property also holds for $J$ and for $F(U)-J$. This together with the fact that the selenoidal property is a linear invariant (i.e.~the set of all states that satisfy the discrete selenoidal constraint forms a vector space which we denote by $\mathcal{V}$) allows us to apply Theorem 2 in \cite{einkemmer1408} which implies that any exponential Runge--Kutta method satisfies the selenoidal property. The only issue here is that the cited theorem assumes that the computation of $\varphi_k(J)v$ preserves the linear constraint in question. However, since we use a Krylov approach the term $\varphi_k(J)v$ is approximated by $p(J)v$, where $p$ is a polynomial. Since we already know that $\mathcal{V}$ is an invariant subspace with respect to $J$, the desired result follows immediately.
%In order to show this result is is sufficient to consider the following ordinary differential equation
%\[ U^{\prime} = F(U), \]
%where $U$ lies in a finite dimensional vector space (i.e.~the space discretization has already been performed by using the classic centered stencil as outlined above). Then it is well known (see, for example, \citep{reynolds2010}) that $F(U)$ satisfies the discrete selenoidal condition
%\[ \frac{B_{i+1,j}(t)-B_{i-1,j}(t)}{2h} + \frac{B_{i,j+1}(t)-B_{i,j-1}(t)}{2h} = 0 \]
%if $U$ does. This implies at once that the same is true forthis ordinary differential equation preserves the discrete selenoidal condition (for all $t\geq 0$)
%if the same is true for the initial value. %In addition, the numerical
%simulations conducted in the next section suggest that the divergence free condition is preserved if an exponential
%integrator is used for the time discretization.
}

\section{The reconnection problem\label{sec:The-reconnection-problem}}

The examples in this and the next section are drawn from \citep{reynolds2006}
and \citep{reynolds2010}, respectively. We start with a reconnection
problem for which the initial value of the magnetic field is given
by

\[
B_{0}(x,y,z)=\left[\begin{array}{c}
\tanh(2y)-\psi_{0}k_{y}\cos(k_{x}x)\sin(k_{y}y)\\
\psi_{0}k_{x}\sin(k_{x}x)\cos(k_{y}y)\\
0
\end{array}\right],
\]
where as in \citep{reynolds2006} we have chosen $k_{x}=\pi/x_{r}$,
$k_{y}=\pi/(2y_{r})$, $\psi_{0}=0.1$, and the computational domain
is given by $[-x_{r},x_{r}]\times[-y_{r},y_{r}]$ for $x_{r}=12.8$
and $y_{r}=\wasblue{6.4}$. This implies that the magnetic field reverses direction
from pointing along $\boldsymbol{e}_{x}$ to pointing along $-\boldsymbol{e}_{x}$
abruptly at $y=0$. Furthermore, \wasblue{at time $t=0$,} we impose a density that is given
by
\[
\rho=1.2-\tanh^{2}(2y)
\]
and a pressure that is proportional to the density; to be more precise
$p=0.5\rho$ (from which the energy is determined). A vanishing velocity
in both space directions is prescribed.

\begin{figure}
	\begin{center}
	\includegraphics[width=6.5cm]{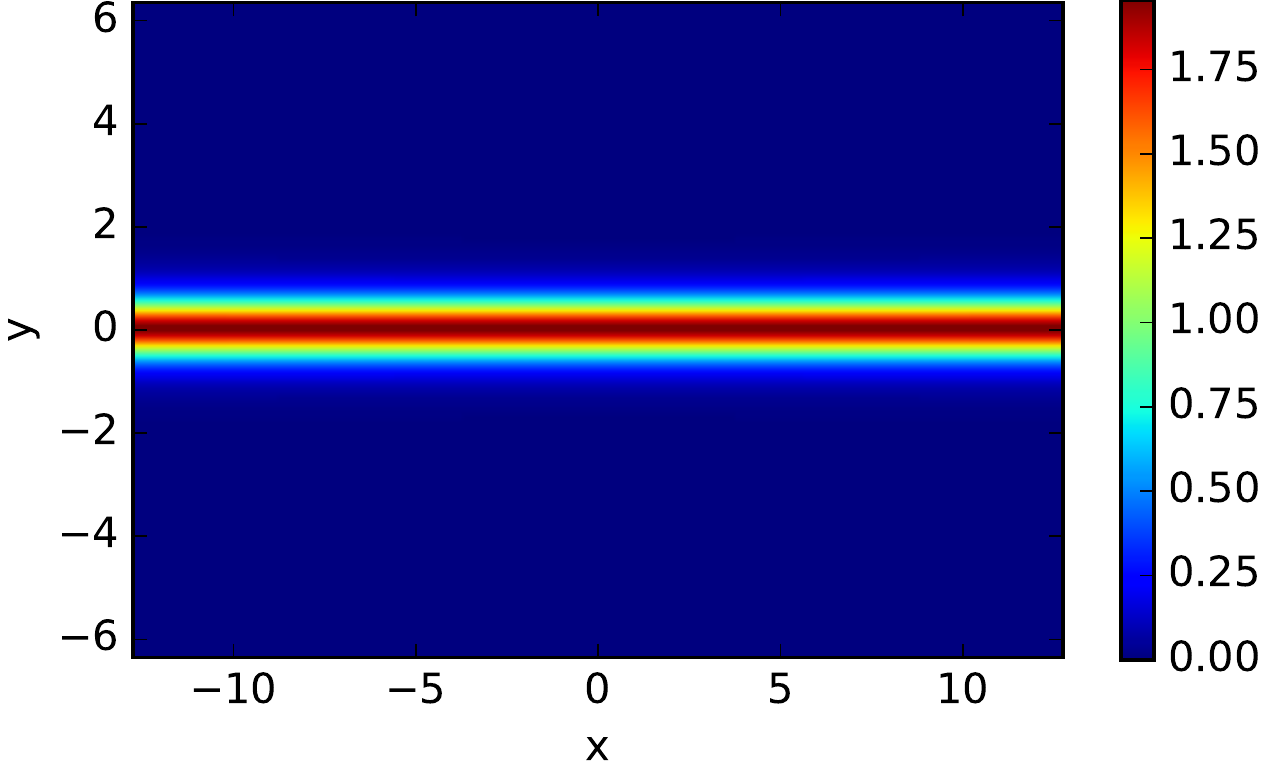}
	\includegraphics[width=6.5cm]{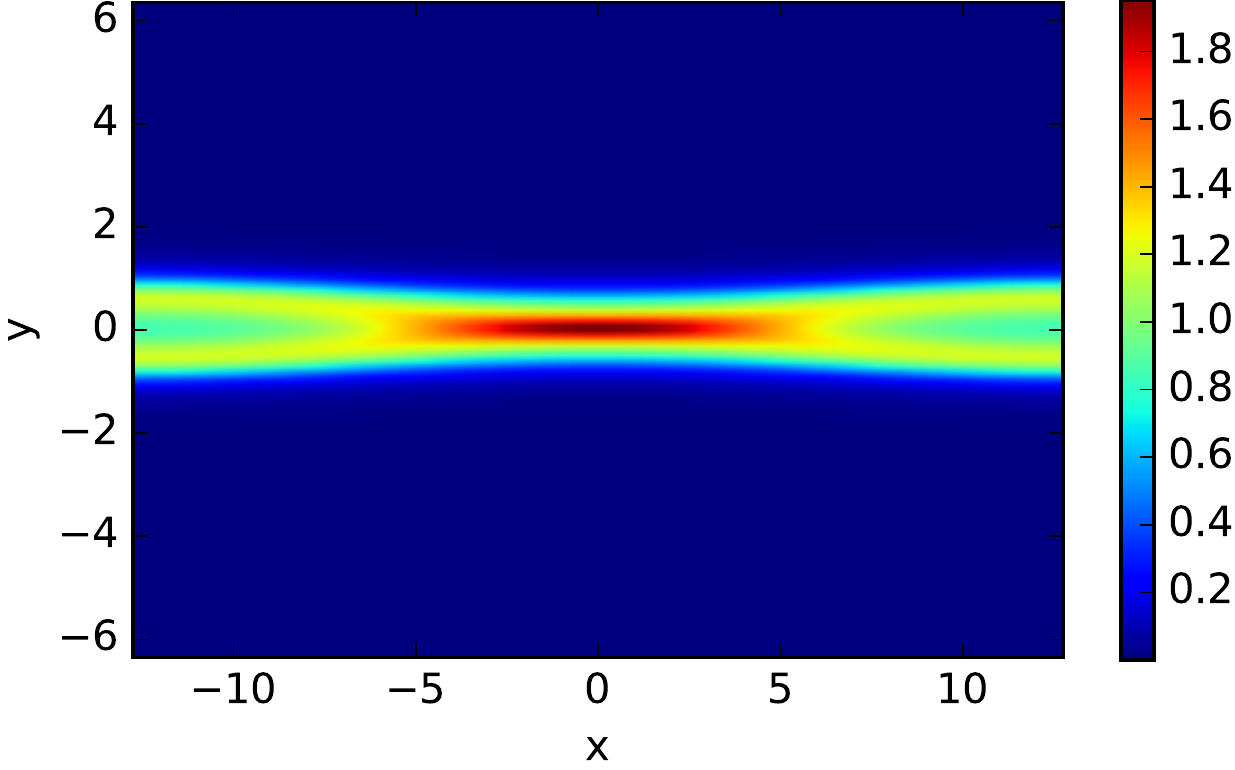}

	\includegraphics[width=6.5cm]{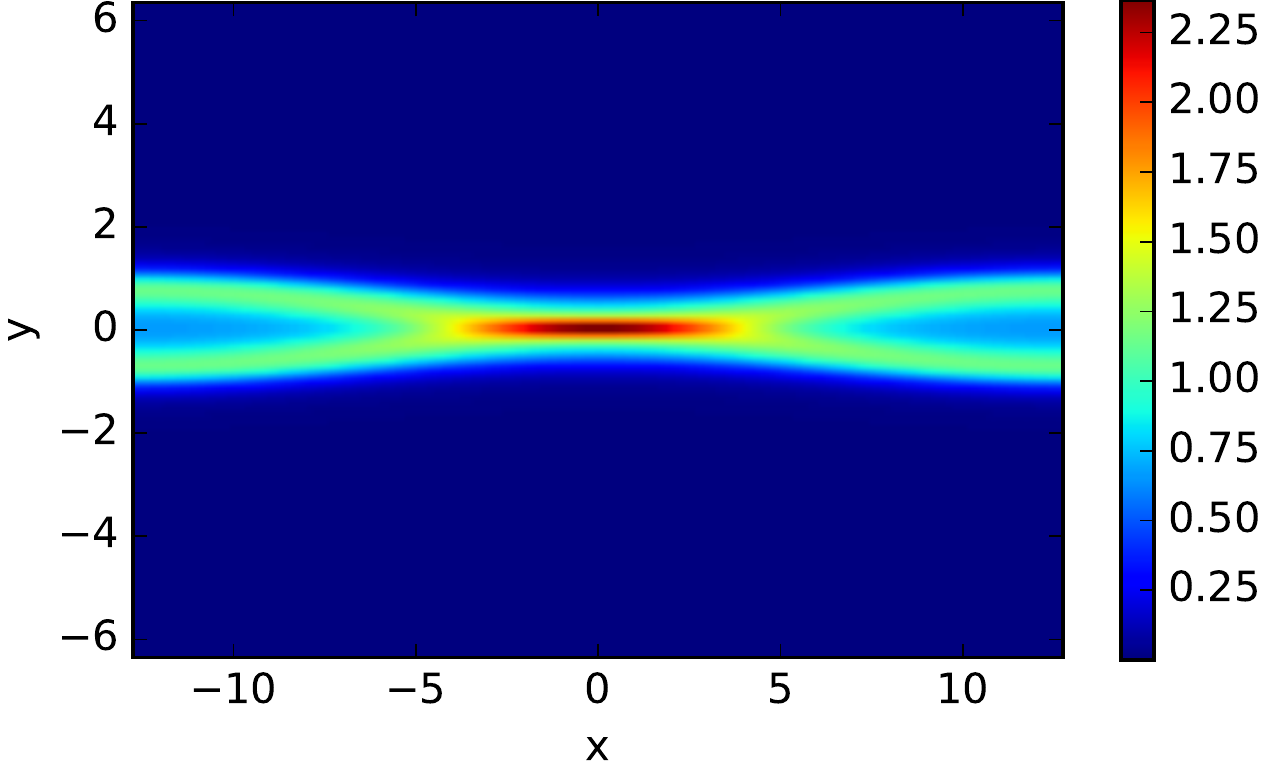}
	\includegraphics[width=6.5cm]{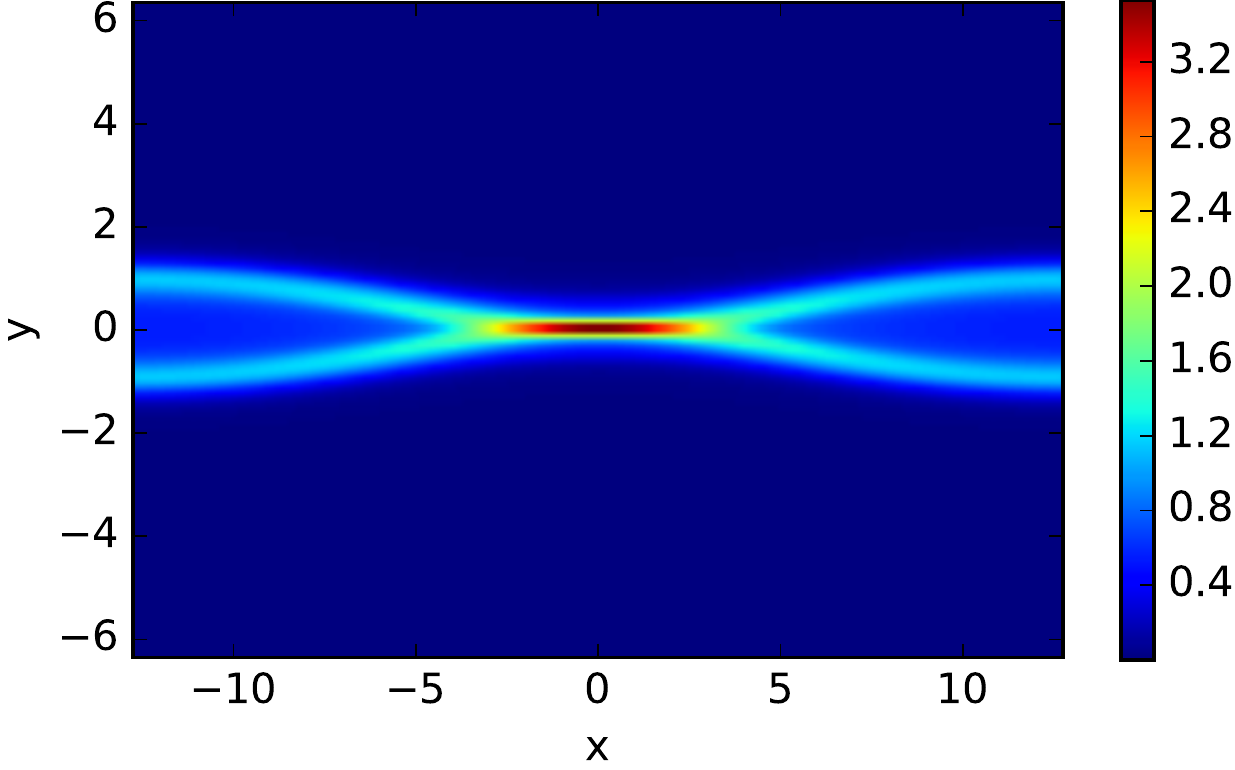}
	
	\includegraphics[width=13cm]{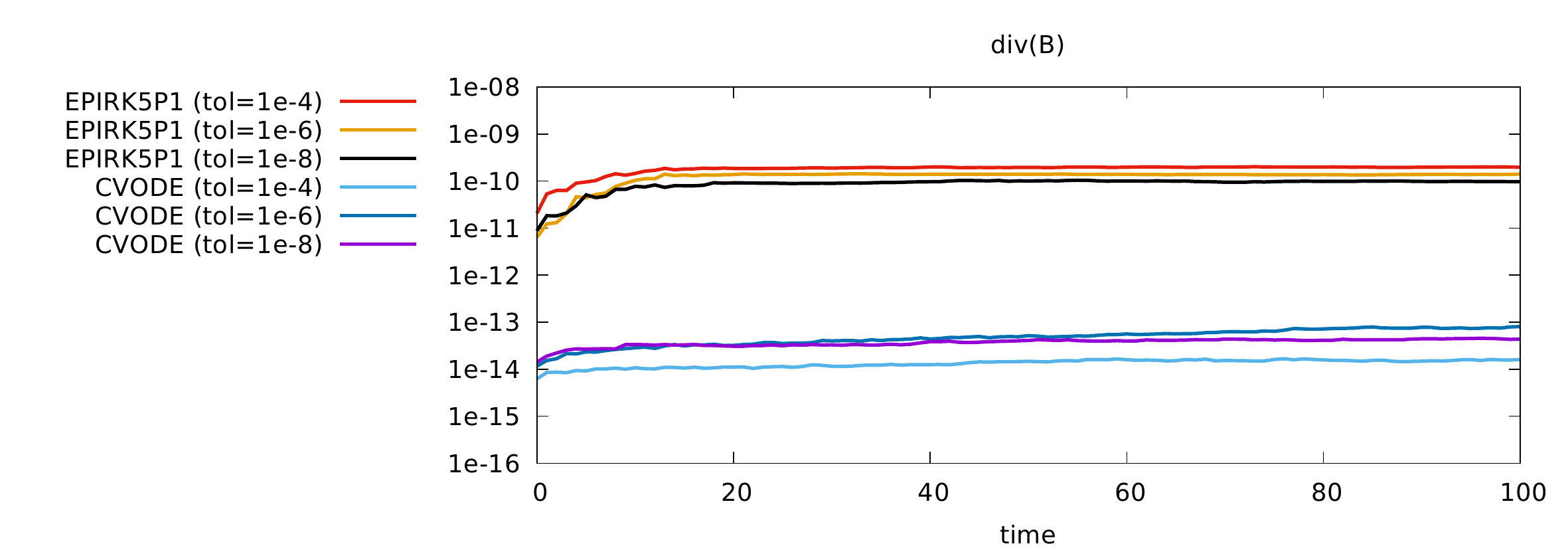}
	\end{center}
	\caption{The current \wasblue{$J=\nabla\times B$} for the reconnection problem is shown at times $t=0$ (top-left), $t=100$ (top-right), $t=150$ (\wasblue{middle}-left), and $t=200$ (\refone{middle}-right). We employ $256$ grid points in the $x$-direction and $128$ grid points in the $y$-direction. The dimensionless parameters are chosen as follows: $\mu=10^{-2}$, $\eta=10^{-3}$, and $\kappa=10^{-2}$. The absolute and relative tolerance for EPIRK5P1 is set to $10^{-6}$. \wasred{On the bottom graph $\nabla\cdot B$ is shown as a function of time for a variety of tolerances.} \label{fig:recon-snapshots}}
\end{figure}

The time evolution for the current is shown in Figure \ref{fig:recon-snapshots}. In all the simulations conducted the relative tolerance, the absolute tolerance, and the tolerance for the Krylov iteration are equal. We have investigated the effect of varying the tolerance for the Krylov iteration and found no significant difference in accuracy or run time as long as we choose a tolerance for the Krylov iteration that is at least as small as the tolerance for the time integration. \wasred{In addition, Figure \ref{fig:recon-snapshots} shows that the value of $\nabla \cdot B$ is (almost) independent of the tolerance chosen and therefore that the error made is only due to numerical round-off.}

\begin{figure}
\noindent \begin{centering}
\includegraphics[width=10cm]{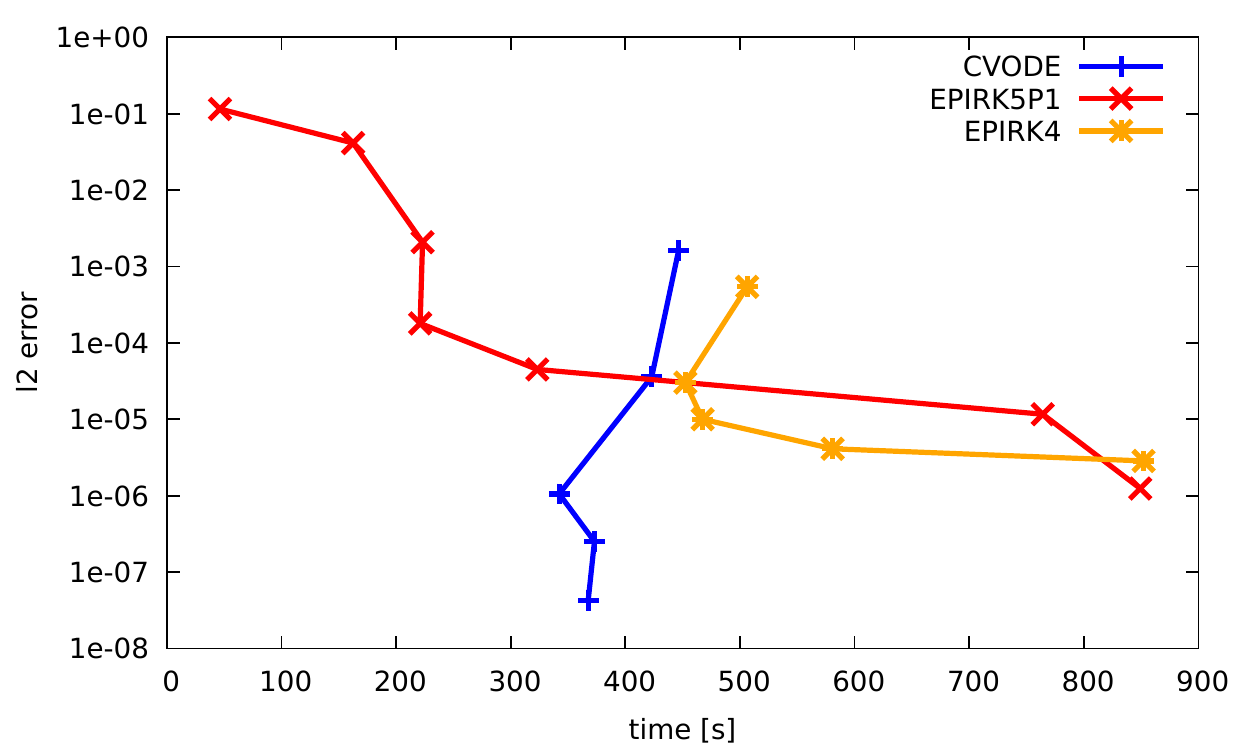}
\par\end{centering}

\protect\caption{The maximum error of the components of the numerical solution (\wasblue{measured
using the $l^{2}$ norm in space at the final time}) is shown as a function of the run time. The
reference solution for the CVODE/EpiRK is computed using the EpiRK/CVODE
method with a tolerance of $10^{-14}$. The simulation is conducted
up to a final time of $T=100$. We have used $256$ grid points in
the $x$-direction and $128$ grid points in the $y$-direction (this
corresponds to the configuration investigated in \citep{reynolds2006}).
The viscosity is given by $\mu=5\cdot10^{-2}$, the resistivity by
	$\eta=5\cdot10^{-3}$, and the thermal conductivity by $\kappa=4\cdot10^{-2}$.
	CVODE tolerance: $10^{-4}$, $10^{-5}$, $10^{-7}$, $10^{-8}$, $10^{-9}$. EPIRK5P1 tolerance: $10^{-1}$, $2\cdot 10^{-2}$, $10^{-2}$, $10^{-3}$, $10^{-5}$, $10^{-7}$, $10^{-9}$. EPIRK4 step size: $2$, $1$, $4\cdot 10^{-1}$, $2\cdot 10^{-1}$, \wasred{$10^{-1}$}.
\label{fig:reconnection-256x128}}
\end{figure}

First, we consider the configuration that is investigated in \citep{reynolds2006};
that is, we employ $256$ grid points in the $x$-direction and $128$
grid points in the $y$-direction. The results obtained are shown in Figure \ref{fig:reconnection-256x128}.
We observe that the performance of the EPIRK5P1 method is superior to CVODE if the desired accuracy is equal or less than $10^{-4}$. In the low precision regime that is often of interest in applications, EPIRK5P1 outperforms CVODE by almost a factor of $2$. However, for more stringent tolerances the CVODE implementation is significantly more efficient. This is due to the fact that the run time for the CVODE method is almost independent of the accuracy achieved. What is perhaps even more interesting is that for a specified tolerance of $10^{-4}$ or above the CVODE method does not converge (in the situation described a final step size on the order of $10^{-16}$ is reported by the application). Let us also note that the constant time step method EPIRK4 is not competitive for this problem.

Next, we investigate the effect of increasing the number of grid points. The results obtained with $512$ grid points in the $x$-direction and $256$ grid points in the $y$-direction are shown Figure \ref{fig:reconnection-512x256}. In this case almost identical conclusions can be drawn. Once again the performance of EPIRK5P1 is superior to CVODE library if the desired accuracy is equal or less than $10^{-4}$. Furthermore, the constant step size method EPIRK4 is not competitive except for very stringent tolerance requirements.

\noindent \begin{center}
\begin{figure}
\noindent \begin{centering}
\includegraphics[width=10cm]{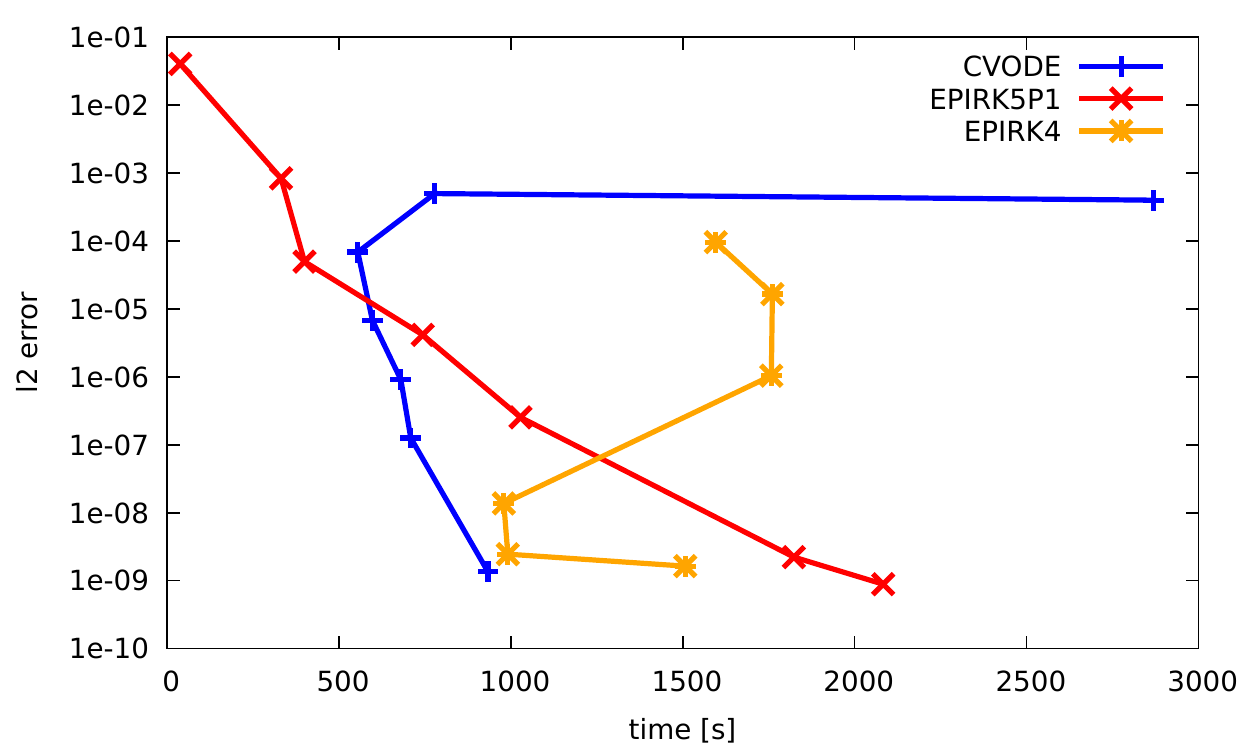}
\par\end{centering}

\protect\caption{The maximum error of the components of the numerical solution
(\wasblue{measured using the $l^{2}$ norm in space at the final time}) is shown as a function of the run time. The
reference solution for the CVODE/EpiRK is computed using the EpiRK/CVODE
method with a tolerance of $10^{-11}$. The simulation is conducted
up to a final time of $T=20$. We have used $512$ grid points in
the $x$-direction and $256$ grid points in the $y$-direction. The
viscosity is given by $\mu=5\cdot10^{-2}$, the resistivity by $\eta=5\cdot10^{-3}$,
and the thermal conductivity by $\kappa=4\cdot10^{-2}$.
	CVODE tolerance: $10^{-3}$, $10^{-4}$, \dots, $10^{-9}$. EPIRK5P1 tolerance: $10^{-1}$, $10^{-2}$, $10^{-3}$, $10^{-4}$, $10^{-5}$, $10^{-8}$, $10^{-9}$. EPIRK4 step size: $4$, $2$, $1$, $4\cdot10^{-1}$, $2\cdot 10^{-1}$, $10^{-1}$.\label{fig:reconnection-512x256}}
\end{figure}

\par\end{center}

Up to this point we have only considered relatively large values for the viscosity, resistivity, and the thermal conductivity. Thus, we proceed by decreasing these dimensionless quantities by a factor of five. The corresponding results are shown in Figure \ref{fig:reconnection-eta-1}. In this case the EPIRK5P1 implementation outperforms the CVODE implementation for all the tolerances studied here (although the difference in performance between the two methods in the low precision regime is smaller than was the case for the previous two examples). For accuracies between $10^{-6}$ and $10^{-8}$ the constant time step method EPIRK4 manages to outperform both the EPIRK5P1 as well as CVODE implementation.

We observe that for larger Reynolds and Lundquist numbers the EpiRK
method shows better or identical performance for the range of tolerances
studied here. In addition, the EpiRK method shows a clear advantage
over the CVODE implementation if a \refone{precision equal to or lower than  $10^{-5}$ is required}.

\noindent \begin{center}
\begin{figure}
\noindent \begin{centering}
\includegraphics[width=10cm]{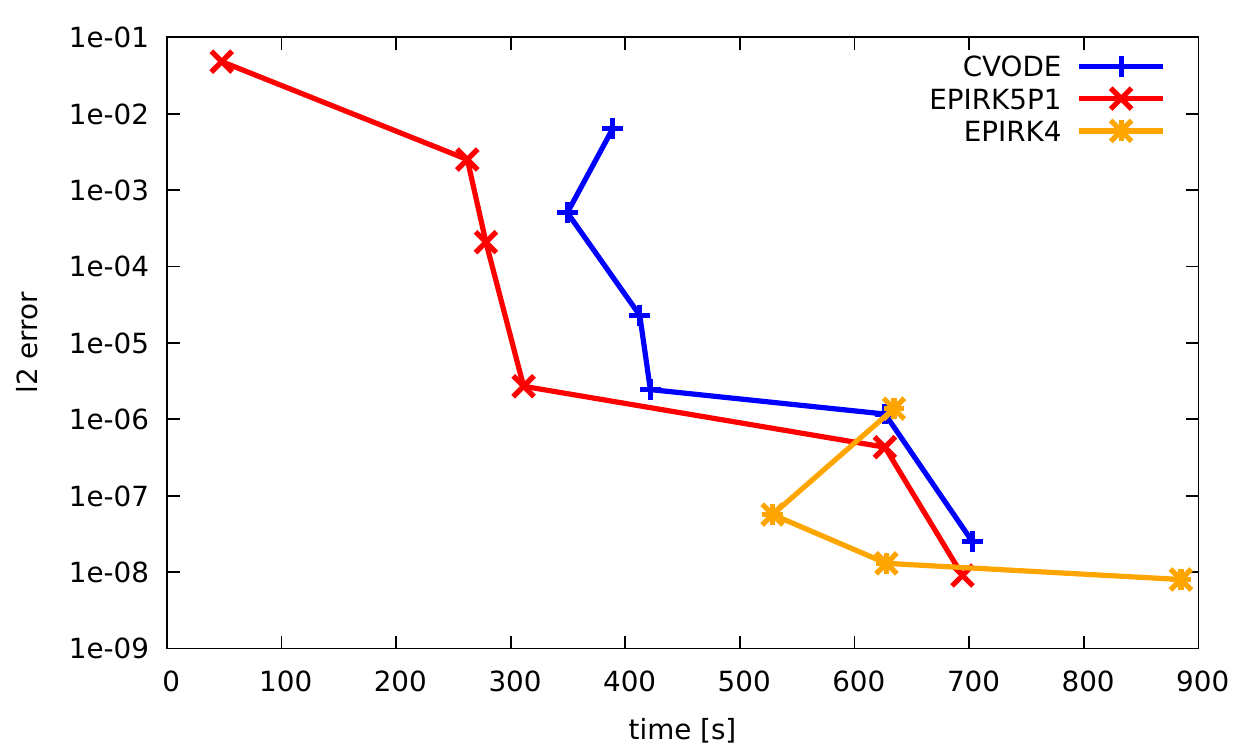}
\par\end{centering}

\protect\caption{The maximum error of the components of the numerical solution
(\wasblue{measured using the $l^{2}$ norm in space at the final time})
is shown as a function of the run time. The
reference solution for the CVODE/EpiRK is computed using the EpiRK/CVODE
method with a tolerance of $10^{-14}$. The simulation is conducted
up to a final time $T=100$. We have used $256$ grid points in the
$x$-direction and $128$ grid points in the $y$-direction. The viscosity
is given by $\mu=10^{-2}$, the resistivity by $\eta=10^{-3}$, and
the thermal conductivity by $\kappa=10^{-2}$.
	CVODE tolerance: $10^{-4}$, $10^{-5}$, \dots, $10^{-9}$. EPIRK5P1 tolerance: $10^{-1}$, $10^{-2}$, $10^{-3}$, $10^{-5}$, $10^{-7}$, $10^{-8}$. EPIRK4 step size: $2$, $1$, $2\cdot10^{-1}$, $10^{-1}$.
    Note that CVODE does not converge for tolerances \refone{that are less tight than $10^{-4}$}.
	\label{fig:reconnection-eta-1}}
\end{figure}

\par\end{center}

In summary, it can be said that for the reconnection problem the EPIRK5P1 implementation shows superior performance in the \refone{loose} and medium tolerance range that is of interest for the majority of practical applications. However, for some of the examples considered here the CVODE implementation has a significant advantage for more stringent tolerance requirements.

To conclude this section let us investigate the strong scaling (i.e.~the problem size is kept
fixed while the number of cores is increased) for both our time integration routines and the CVODE implementation. This allows us to ascertain how well the current implementation is parallelized. In the ideal case we would observe a linear decrease in the run time as a function of the number of cores.

The strong scaling results on the VSC-2 are shown in Figure \ref{fig:scaling} and are almost identical for both implementations. We conclude that the algorithms scale well up to approximately $256$ cores. No further gain in performance can be observed for $1024$ or more cores. This is to be expected as in the case of strong scaling the amount of work available to a single core eventually decreases sufficiently such that communication between he different MPI processes limits the performance of the application. At this point no further increase in performance can be expected.

\noindent \begin{center}
\begin{figure}
\noindent \begin{centering}
\includegraphics[width=10cm]{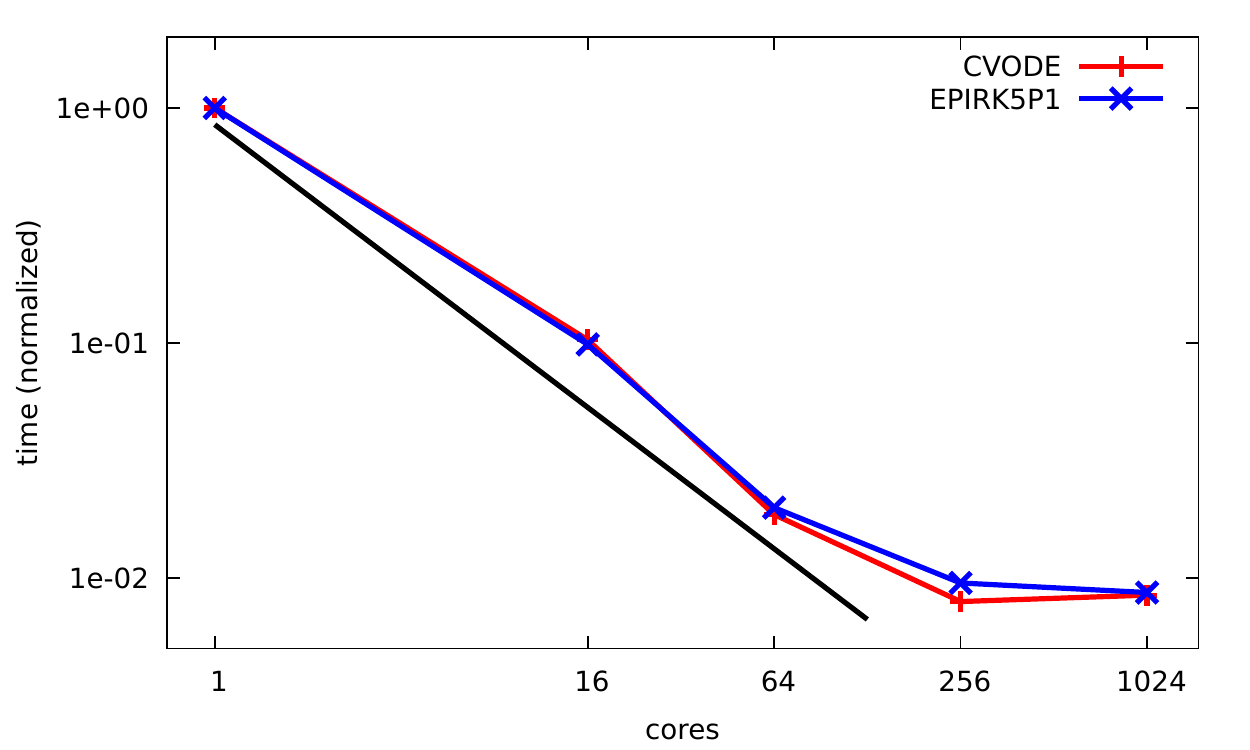}
\par\end{centering}

\protect\caption{The run time (normalized to the single core run time for both the
CVODE and EpiRK implementation) is shown as a function of the number
of cores. The simulation is conducted up to a final time $T=1$. We
have used $512$ grid points in both the $x$- and $y$-direction.
The viscosity is given by $\mu=5\cdot10^{-2}$, the resisitivity by
	$\eta=5\cdot10^{-3}$, and the thermal conductivity by $\kappa=4\cdot10^{-2}$. These performance measurements were done on the VSC-2 (Vienna Scientific Cluster 2; see \protect\url{http://vsc.ac.at/systems/vsc-2/}) which uses two 8 core AMD Opteron Magny Cours 6132HE and $32$ GB of DDR3 memory per node. The nodes are connected by Infiniband QDR and the Intel MPI library (version 4.1) is employed. \wasblue{The line of ideal scaling is shown in black.}
\label{fig:scaling}}
\end{figure}

\par\end{center}

\section{The Kelvin--Helmholtz instability\label{sec:The-Kelvin-Helmholtz-instability}}

As a second example we consider the Kelvin--Helmholtz (KH) instability.
The KH instability is triggered by superimposing the following perturbation
\[
\epsilon_{x}\cos\left(\tfrac{2\pi\omega_{x}}{L_{x}}x\right)+\epsilon_{y}\sin\left(\frac{\pi(2\omega_{y}-1)}{L_{y}}y\right)
\]
 in the $x$-direction on the velocity field
\[
\boldsymbol{v}=\left[\begin{array}{c}
v_{0}\tanh(\tfrac{y}{\lambda})\\
0\\
0
\end{array}\right].
\]
 The density is initialized to unity and a uniform pressure is chosen.
The magnetic field is initialized to be uniform in both the $x$-
and $z$-direction with strength $B_{x}$ and $B_{z}$, respectively,
and is assumed to vanish in the $y$-direction. All the parameters
used to determine the numerical value of the initial value are listed
in Table \ref{tab:kh-parameters}.

\begin{figure}
	\begin{center}
		\includegraphics[width=6.5cm]{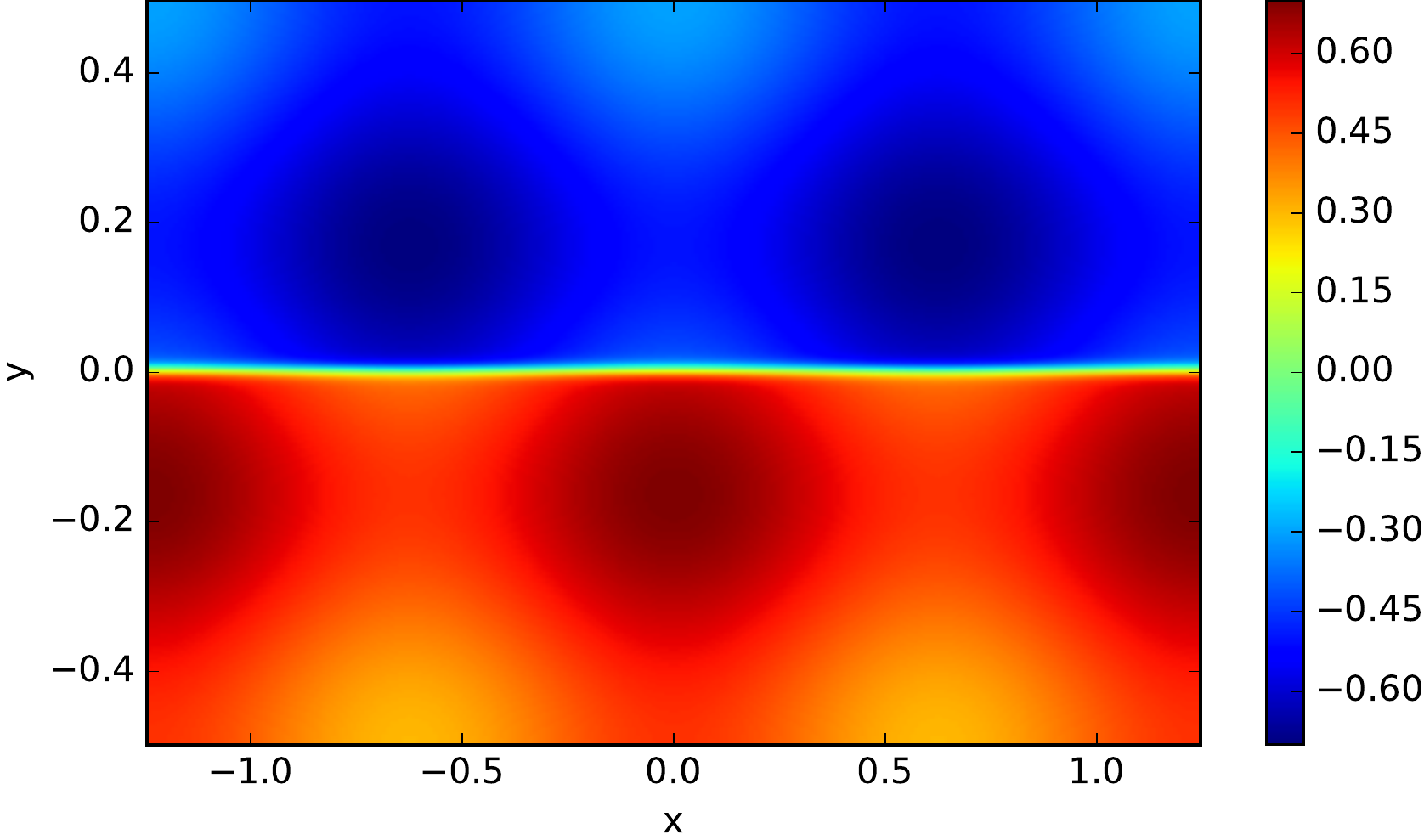}
		\includegraphics[width=6.5cm]{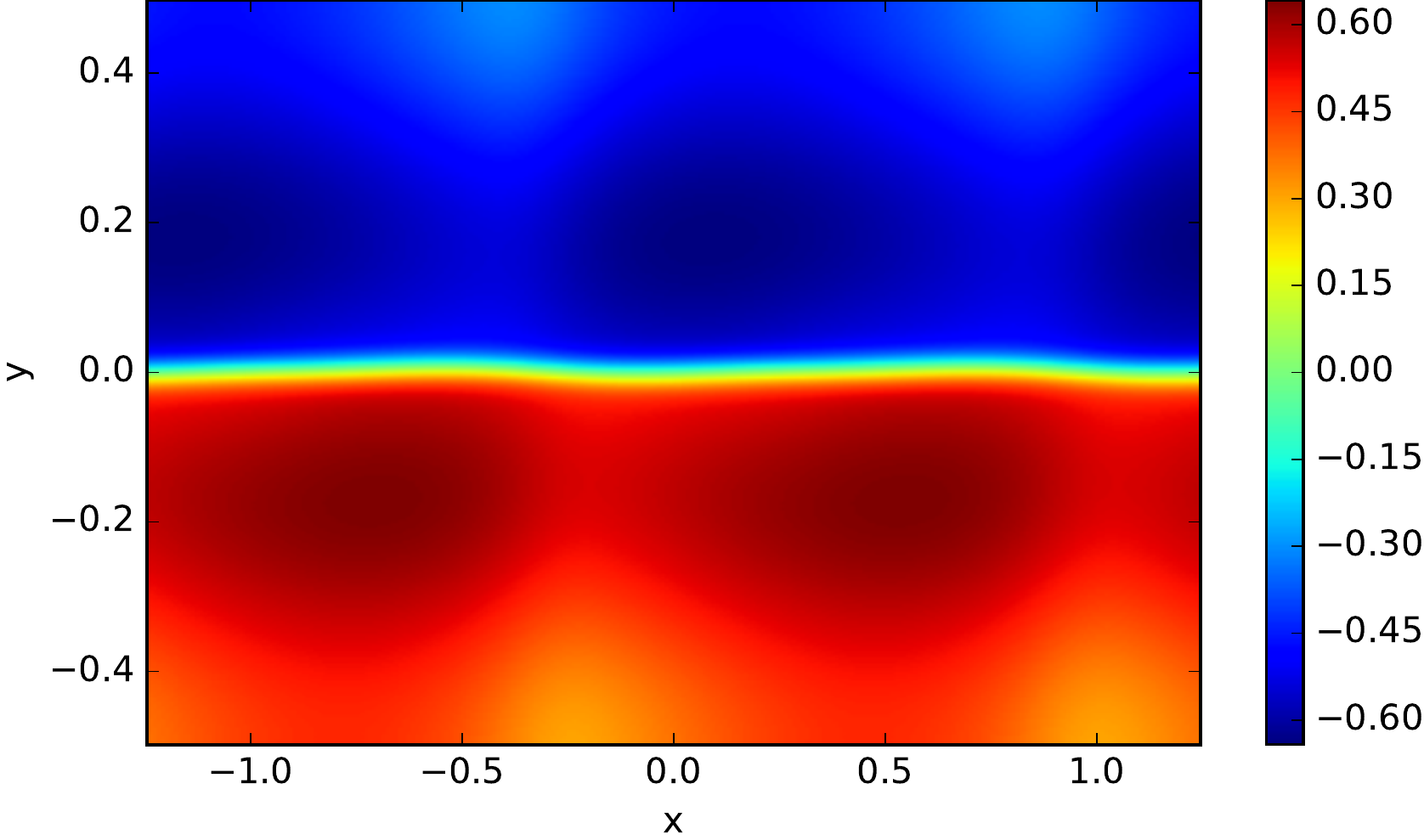}
		
		\includegraphics[width=6.5cm]{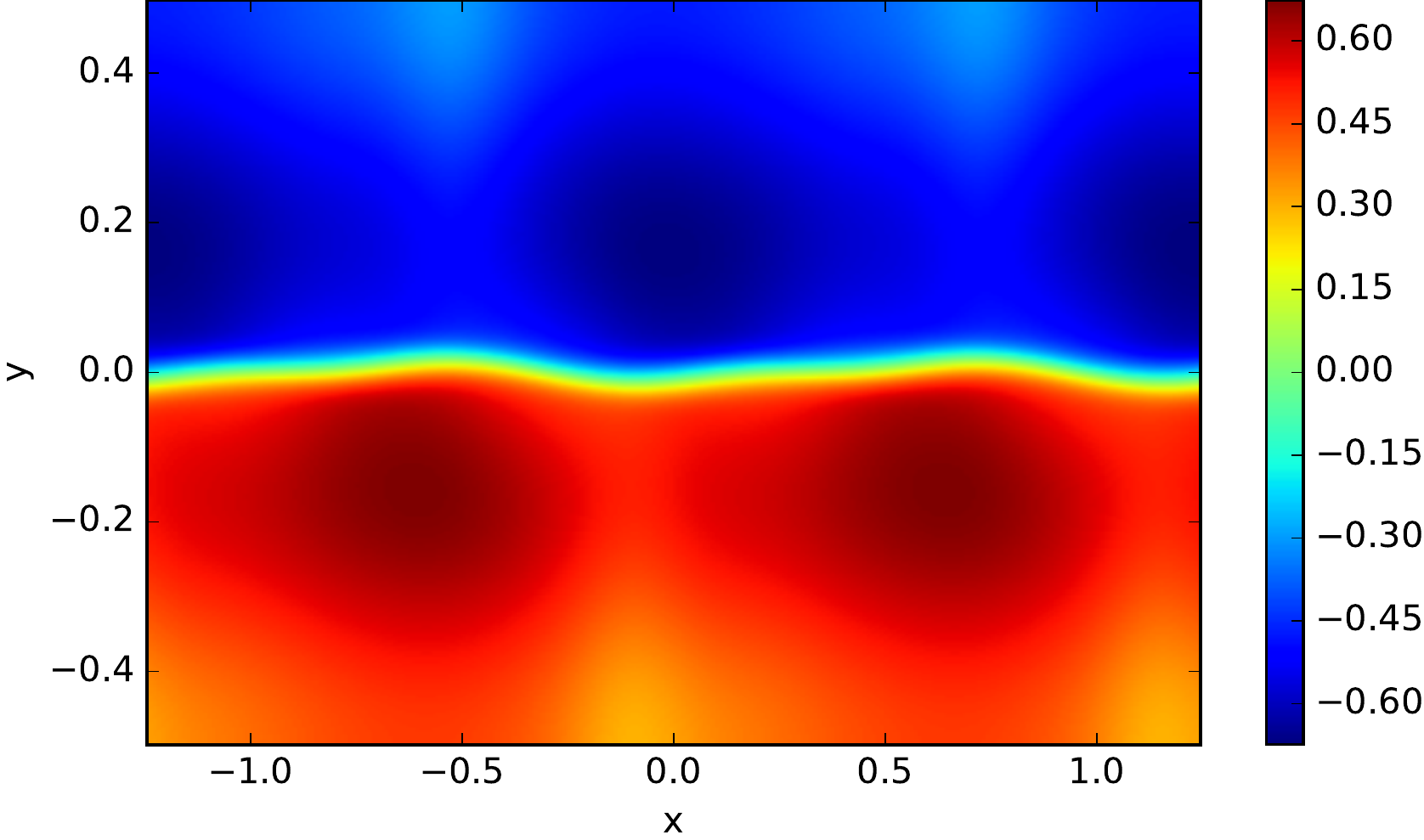}
		\includegraphics[width=6.5cm]{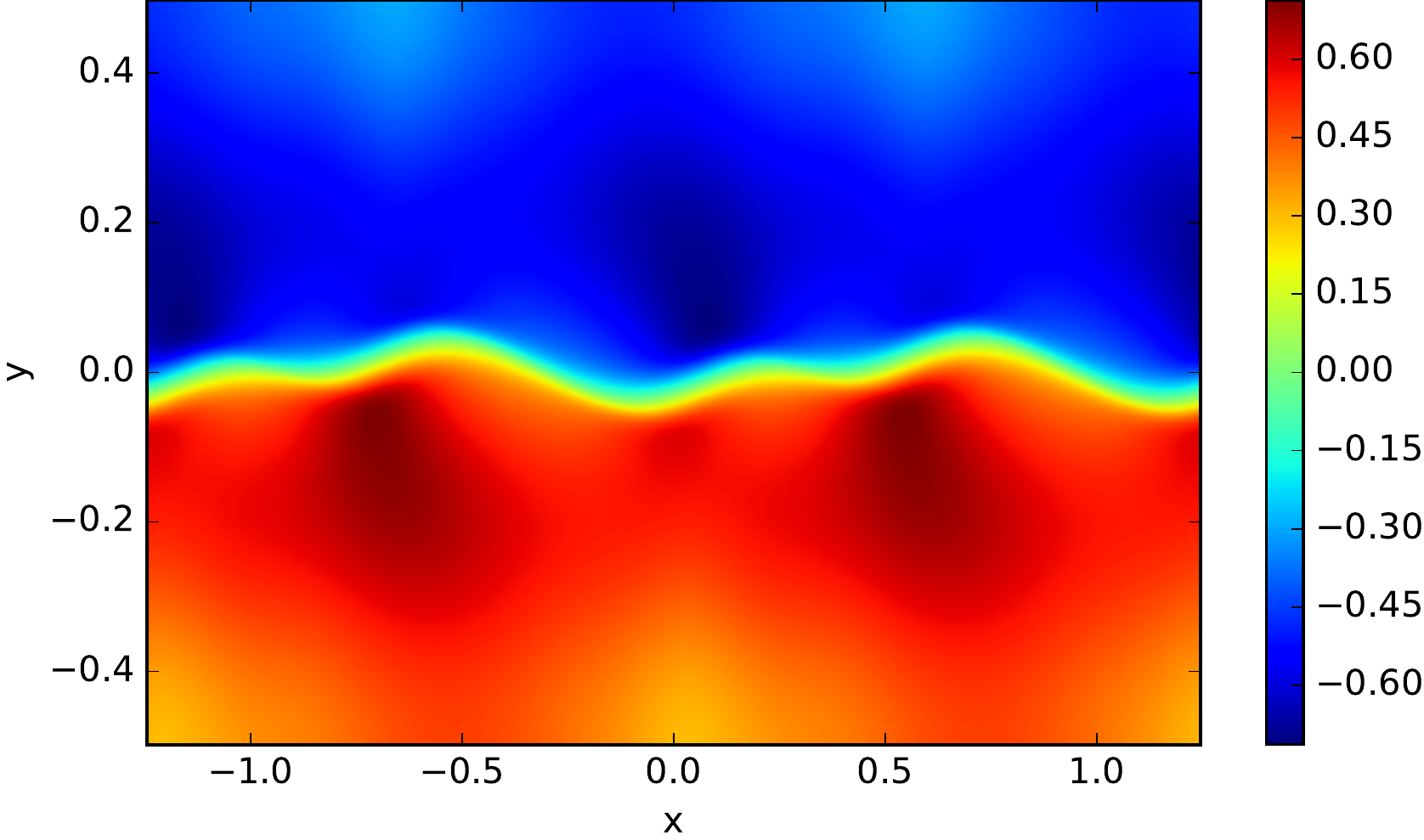}
	\end{center}
	\caption{For the Kelvin--Helmholtz instability the velocity in the $x$-direction is shown at time $t=0$ (top-left), $t=1.5$ (top-right), $t=2$ (bottom-left), and at $t=2.5$ (bottom-right). In both space directions $256$
  grid points are employed. The dimensionless parameters are chosen as follows: $\mu=10^{-4}$, $\eta=10^{-4}$, and $\kappa=10^{-4}$. The absolute and relative tolerance for EPIRK5P1 is set to $10^{-6}$.\label{fig:kh-snapshots}}
\end{figure}

The time evolution of the numerical solution is shown in Figure \ref{fig:kh-snapshots}. In all the simulations conducted the relative tolerance, the absolute tolerance, and the tolerance for the Krylov iteration are equal. We have investigated the effect of varying the tolerance for the Krylov iteration and found no significant difference in accuracy or run time as long as we choose a tolerance for the Krylov iteration that is at least as small as the tolerance for the time integration.

\begin{table}
\noindent \begin{centering}
\begin{tabular}{cc}
parameter & value\tabularnewline
\hline
$\epsilon_{x},\epsilon_{y}$ & $0.1$\tabularnewline
$v_{0}$ & $0.5$\tabularnewline
$\omega_{x},\omega_{y}$ & $2$\tabularnewline
pressure $p$ & $0.25$\tabularnewline
$B_{z}$ & $10$\tabularnewline
$B_{x}$ & $0.1$\tabularnewline
\end{tabular}
\par\end{centering}

\protect\caption{Parameters for the initial value of the Kelvin--Helmholtz instability.
 \label{tab:kh-parameters}}
\end{table}

\begin{figure}
\noindent \begin{centering}
\includegraphics[width=10cm]{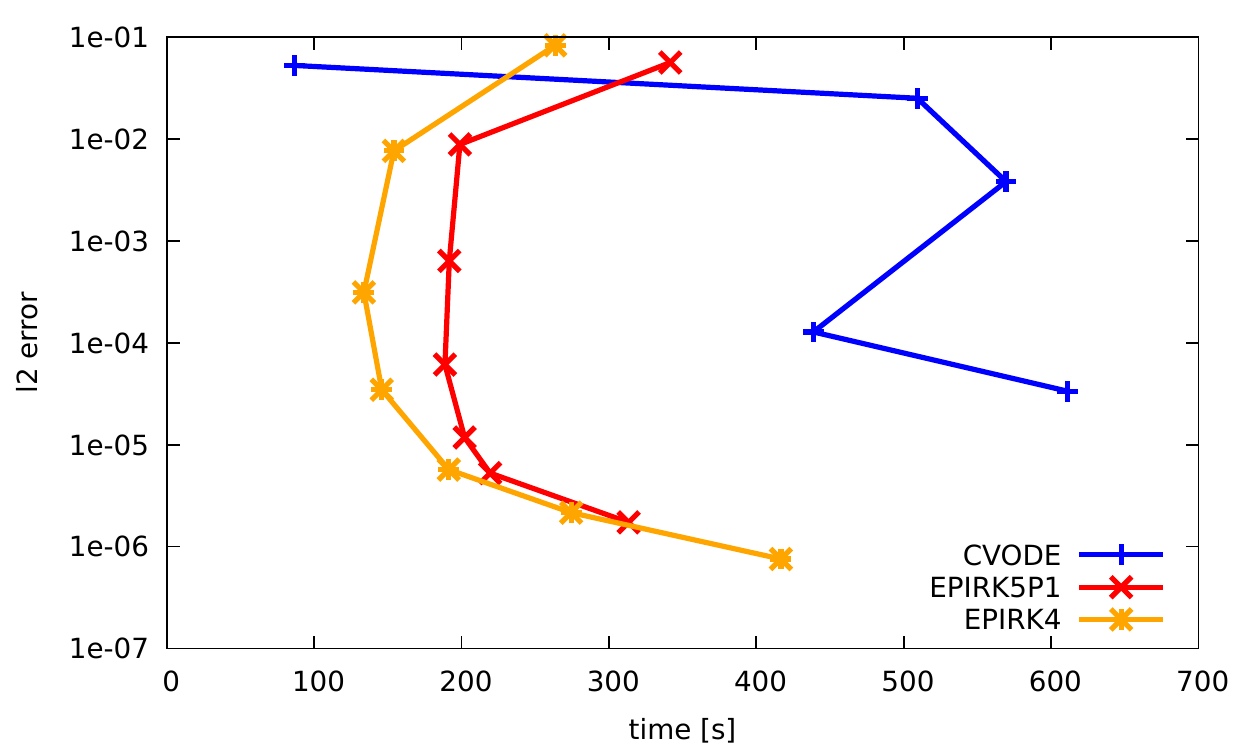}
\par\end{centering}

\noindent \begin{centering}
\includegraphics[width=10cm]{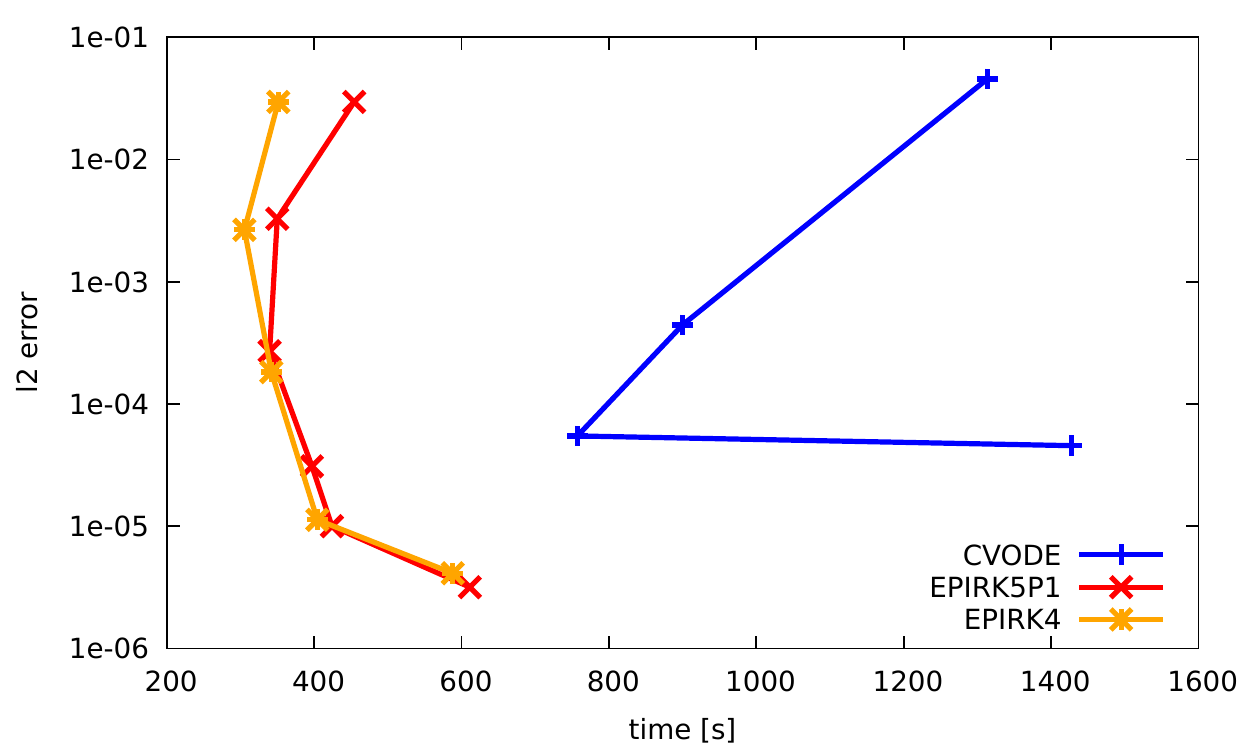}
\par\end{centering}

\protect\caption{The maximum error of the components of the numerical solution
(\wasblue{measured using the $l^{2}$ norm in space at the final time})
is shown as a function of the run time. The
reference solution for the CVODE/EpiRK is computed using the EpiRK/CVODE
method with a tolerance of $10^{-11}$. The simulation is conducted
up to a final time $T=1$ (top) and $T=2$ (bottom) and $128$ grid
points are used in both the $x$- and $y$-direction. The (dimensionless)
viscosity is given by $\mu=10^{-4}$, the resistivity by $\eta=10^{-4}$,
and the thermal conductivity by $\kappa=10^{-4}$.
	CVODE tolerance ($T=1$): $10^{-3}$, $10^{-6}$, $10^{-7}$, $10^{-8}$, $10^{-9}$.
	CVODE tolerance ($T=2$): $10^{-6}$, $10^{-7}$, $10^{-8}$, $10^{-9}$.
	EPIRK5P1 tolerance ($T=1$): $10^{-3}$, $10^{-4}$, \dots, $10^{-9}$.
	EPIRK5P1 tolerance ($T=2$): $10^{-4}$, $10^{-5}$, \dots, $10^{-9}$.
	EPIRK4 step size ($T=1$): $10^{-1}$, $4\cdot 10^{-2}$, $2\cdot 10^{-2}$, $10^{-2}$, $4\cdot 10^{-3}$, $2\cdot 10^{-3}$, $10^{-3}$.
	EPIRK4 step size ($T=2$): $4\cdot 10^{-2}$, $2\cdot 10^{-2}$, $10^{-2}$, $4\cdot 10^{-3}$, $2\cdot 10^{-3}$.
	\label{fig:kh-128}}
\end{figure}

First, we perform a numerical simulation of the Kelvin--Helmholtz instability for $128$ grid points in both space directions. The results obtained are shown in Figure \ref{fig:kh-128}. Both exponential integrators considered here (EPIRK5P1 and EPIRK4) outperform the CVODE implementation by a significant margin (approximately a factor of two). Furthermore, CVODE is very sensitive to even small perturbations in the chosen tolerance, while the exponential integrators are much more robust in this regard. Thus, a user who would naively (i.e.~without measuring the run time for a range of tolerances) choose a specific tolerance would increase the run time of the application by significantly more than the factor of two stated earlier. What is also interesting to see is that the constant step size method EPIRK4 is equal in performance or outperforms EPIRK5P1 for all tolerances studied here (although the difference between these two methods is at most 30\%).

Now, let us increase the number of grid points to $256$ (in both space directions). The corresponding results are shown in Figure \ref{fig:KH-512}. As before, we observe that both exponential integrators significantly outperform the CVODE implementation. In addition, the constant step size method EPIRK4 outperforms the EPIRK5P1 implementation by up to 30\%. This is also consistent with the results on the coarse grid considered before. Let us emphasize here again that the major advantage of EPIRK4 is that it uses fewer Krylov projections (two vs. three for EPIRK5P1) per time step.  However, at present this advantage over EPIRK5P1 is limited to constant time step implementations.  A third projection is needed to construct automatic time stepping control schemes for EPIRK4.  The variable time stepping version of the integrators is more efficient for EPIRK5P1 since it is a higher order method. For the Kelvin--Helmholtz instability, considered in this section, these differences work in favor of the EPIRK4 method while for the reconnection problem, considered in the previous section, the EPIRK5P1 scheme is superior (in the case of \refone{loose} and medium tolerances).

In summary, both exponential integrators considered in this paper show significantly improved performance compared to CVODE for the Kelvin--Helmholtz instability. Let us also mention that for both the exponential integrators as well as CVODE it is not entirely \wasblue{straightforward} to choose a good tolerance. This is due to the fact that if the tolerance is chosen too \refone{loose}, the run time can actually increase. \wasred{This is due to the deficiency of the step size controller which assumes that the computational cost of each step is independent of the step size (a reasonable assumption for explicit time integrators; however, this is not true for the implicit and exponential methods considered here).} We emphasize, however, that this behavior is significantly more serious for the CVODE implementation.
\noindent \begin{center}
\begin{figure}
\noindent \begin{centering}
\includegraphics[width=10cm]{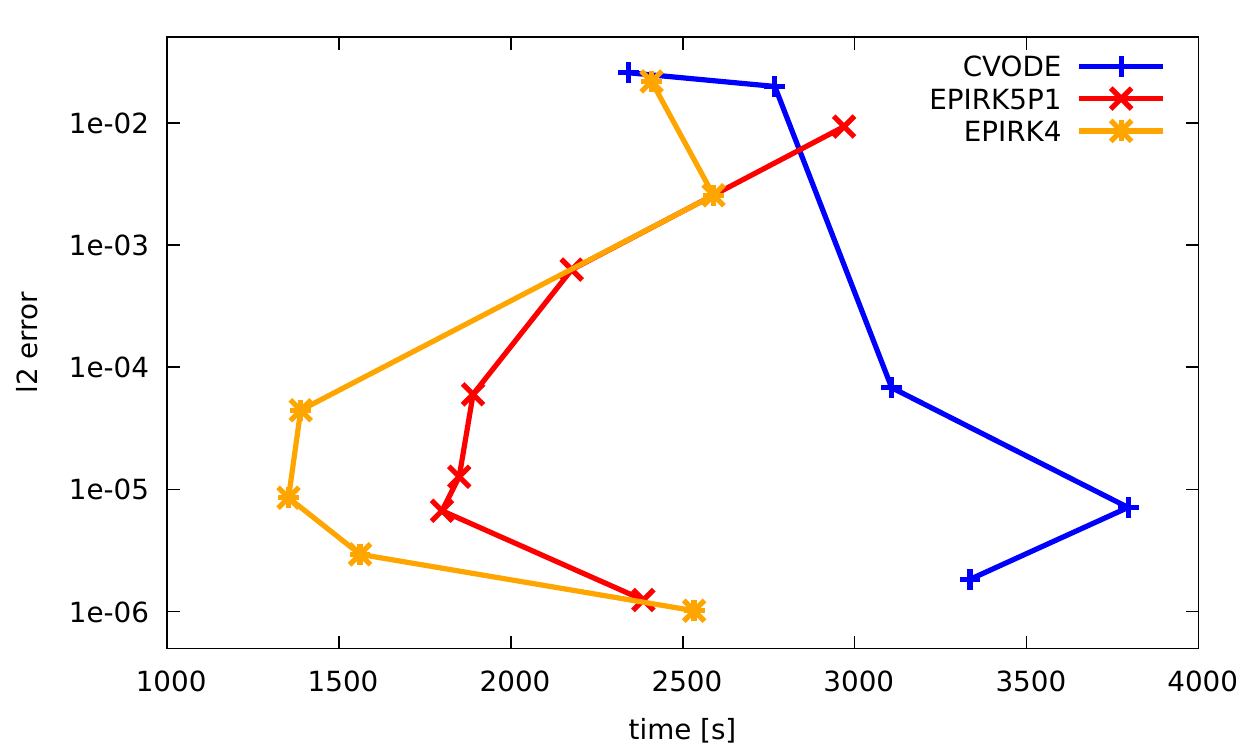}
\par\end{centering}
\protect\caption{The maximum error of the components of the numerical solution
(\wasblue{measured using the $l^{2}$ norm in space at the final time})
is shown as a function of the run time. The
reference solution for the CVODE/EpiRK is computed using the EpiRK/CVODE
method with a tolerance of $10^{-11}$. The simulation is conducted
up to a final time $T=1$ and $256$ grid points are used in both
the $x$- and $y$-direction. The (dimensionless) viscosity is given
by $\mu=10^{-4}$, the resistivity by $\eta=10^{-4}$, and the thermal
conductivity by $\kappa=10^{-4}$.
	CVODE tolerance: $10^{-5}$, $10^{-6}$, \dots, $10^{-9}$.
	EPIRK5P1 tolerance: $10^{-4}$, $10^{-5}$, \dots, $10^{-9}$.
	EPIRK4 step size: $4\cdot 10^{-2}$, $2\cdot 10^{-2}$, $10^{-2}$, $4\cdot 10^{-3}$, $2\cdot 10^{-3}$, $10^{-3}$.
	\label{fig:KH-512}}
\end{figure}
\par\end{center}

\section{Implementation\label{sec:Implementation}}

We have chosen to base our implementation on the $2.5$ dimensional\footnote{All quantities depend on only two spatial variables; however, the
direction of the velocity as well as the magnetic fields are three
dimensional vectors.} MHD code developed by Daniel R. Reynolds et al. A large number of
numerical simulations have been conducted (see e.g. \citep{reynolds2006}
and \citep{reynolds2010}) and the performance of preconditioners
for Newton--Krylov based implicit methods has been investigated (see
e.g. \citep{reynolds2010}). Furthermore, the extension of the code
to three dimensional problems as well as to non-square geometries
(such as a Tokamak geometry) have been investigated in \citep{reynolds2012}.

Our code is implemented in C++ and is designed to easily accommodate different time integrators, space
discretization schemes, as well as different MHD models. The design
considerations and a detailed description of the computer code can
be found in \citep{einkemmercppmhd}.

The second component of our code is the EPIC library (see \citep{loffeld2012})
which implements the fifth order EpiRK5P1 with an adaptive time step control ({\it EpiRK5P1VerticalVar} function in EPIC library),
the fourth order EpiRK4 method ({\it EpiRK4MixedConst} function in EPIC library) with constant time steps and
the adaptive Krylov algorithm to evaluate the $\varphi_{k}$ functions
as described in \citep{niesen2012}. This library is written in the
C++ programming language with MPI, and uses routines
from the BLAS and LAPACK libraries. In the experiments,
the BLAS and LAPACK implementations were supplied by
the Intel Math Kernel Library (Intel MKL).

One detail warrants further discussion: in the original Fortran program
the CVODE interface is used to approximate the Jacobian by a simple
forward difference stencil. In our implementation, however, we provide
a custom function to compute an approximation to the Jacobian. This
is necessary for the EPIC library, which before was only used in the
context of problems where an exact Jacobian is available (as in \citep{loffeld2012}).
Computing an exact Jacobian for the significantly more complex
MHD problem considered in this paper is infeasible. Therefore, we
have to approximate $J(a)v$, i.e. the application of the Jacobian
$J$ at position $a$ to a vector $v$, where both $a$ and $v$ depend
on the specific numerical algorithm as well as the initial values.
We note here that there is a qualitative difference in the norm of
the vector $v$ that depends on the numerical time integration scheme
used. In the EpiRK the norm of $v$ is close to unity whereas in the
BDF method the value is often significantly below $\sqrt{\epsilon}$
(where we use $\epsilon\approx10^{-16}$ to denote machine precision).
Therefore, care is taken to scale vectors to the norm $\sqrt{\epsilon}$
only if the initial norm is above $\sqrt{\epsilon}$ in magnitude.
Then the same implementation can be used for both methods and a difference
in performance due to an internally optimized function to compute
the Jacobian (as is provided by the CVODE library) is precluded.

In all of the computations conducted in this paper, we have used the
infinity norm to scale $v$ to its appropriate size. In general, we
have found that the performance and accuracy of the computation is
not significantly altered if the scaling is performed to some value
that is reasonably close to $\sqrt{\epsilon}$.

\section{Conclusion\label{sec:Conclusion}}

Using several examples of MHD problems, we have shown that exponential
integrators constitute a viable alternative to the more commonly employed
BDF method (as implemented in the CVODE library). In almost all instances
equal or superior performance has been observed for the adaptive and
variable time stepping fifth order EpiRK method for low to medium
accuracy requirements.

Note that while the first version of the CVODE package was released several
decades ago and the package has been optimized and refined over extended period of time, the EPIC package is
the first C++/MPI implementation of the exponential integrators.  As indicated in \cite{loffeld14}
there are a number of additional optimizations possible both from an algorithmic and from a
computer science perspective.  For example, significant performance gains are expected
from improvements in the adaptivity algorithms in the Krylov projections and the derivation of the time
integrators specifically adapted to a given system of equations.  In addition, the
newly introduced implicit-exponential (IMEXP) methods  \cite{rainwater2016} demonstrate that exponential integration
can also be used as a component of the overall integration scheme and can be combined
with the preconditioned implicit integrators as well as time splitting strategies.  It is expected
that further development will lead to significant performance improvements.

Furthermore, we note that the results presented here are somewhat
different from the simpler models considered in \citep{loffeld14}.
In that case order of magnitude speedups are observed for the EpiRK
method as compared to the CVODE implementation. Further investigations
are needed to determine whether the degradation in comparative performance
is due to the specifics of the PDE system, the details of the integrator
or the particulars of the implementation. Also it is not clear to what degree
the approximate computation of the Jacobian contributes to the observed
differences (in all the problems investigated in \citep{loffeld14}
an analytical form of the Jacobian has been used).

\section*{Acknowledgements}

We would like to take the opportunity to thank D. R. Reynolds for providing the
code of the Fortran MHD solver and for the helpful discussion.
Dr. Einkemmer was supported by the Marshall plan scholarship of
the Austrian Marshall plan foundation (http://www.marshallplan.at/)
and by the Fonds zur Förderung der Wissenschaften (FWF) -- project
id: P25346. Dr. Tokman was supported by a grant from the National Science Foundation,
Computational Mathematics Program, under Grant No. 1115978.

The simulations were conducted using the Vienna
Scientific Cluster (VSC).

\bibliographystyle{amsplain}
\bibliography{paper-epic}

\end{document}